\documentclass[a4paper,11pt]{amsart}
\usepackage{mathrsfs}
\usepackage{cases}
\usepackage{amsfonts}
\usepackage{graphicx}
\usepackage{amsmath}
\usepackage{amssymb}
\usepackage[all]{xypic}
\usepackage[all]{xy}
\usepackage{color}
\usepackage{colordvi}
\usepackage{multicol}
\usepackage[colorlinks,citecolor=black,linkcolor=black]{hyperref}
\begin{document}
\input xy
\xyoption{all}

\newcommand{\ind}{\operatorname{inj.dim}\nolimits}
\newcommand{\id}{\operatorname{id}\nolimits}
\newcommand{\Mod}{\operatorname{Mod}\nolimits}
\newcommand{\End}{\operatorname{End}\nolimits}
\newcommand{\Ext}{\operatorname{Ext}\nolimits}
\newcommand{\Hom}{\operatorname{Hom}\nolimits}
\newcommand{\aut}{\operatorname{Aut}\nolimits}
\newcommand{\Ker}{{\operatorname{Ker}\nolimits}}
\newcommand{\Iso}{\operatorname{Iso}\nolimits}
\newcommand{\Coker}{\operatorname{Coker}\nolimits}
\renewcommand{\dim}{\operatorname{dim}\nolimits}
\newcommand{\Cone}{{\operatorname{Cone}\nolimits}}
\renewcommand{\Im}{\operatorname{Im}\nolimits}
\newcommand{\cc}{{\mathcal C}}
\newcommand{\ce}{{\mathcal E}}
\newcommand{\cs}{{\mathcal S}}
\newcommand{\cf}{{\mathcal F}}
\newcommand{\cx}{{\mathcal X}}
\newcommand{\cy}{{\mathcal Y}}
\newcommand{\cl}{{\mathcal L}}
\newcommand{\ct}{{\mathcal T}}
\newcommand{\cu}{{\mathcal U}}
\newcommand{\cm}{{\mathcal M}}
\newcommand{\cv}{{\mathcal V}}
\newcommand{\ch}{{\mathcal H}}
\newcommand{\ca}{{\mathcal A}}
\newcommand{\mcr}{{\mathcal R}}
\newcommand{\cb}{{\mathcal B}}
\newcommand{\ci}{{\mathcal I}}
\newcommand{\cj}{{\mathcal J}}
\newcommand{\cp}{{\mathcal P}}
\newcommand{\cg}{{\mathcal G}}
\newcommand{\cw}{{\mathcal W}}
\newcommand{\co}{{\mathcal O}}
\newcommand{\cd}{{\mathcal D}}
\newcommand{\cn}{{\mathcal N}}
\newcommand{\ck}{{\mathcal K}}
\newcommand{\calr}{{\mathcal R}}
\newcommand{\ol}{\overline}
\newcommand{\ul}{\underline}
\newcommand{\cz}{{\mathcal Z}}
\newcommand{\st}{[1]}
\newcommand{\ow}{\widetilde}
\renewcommand{\P}{\mathbf{P}}
\newcommand{\pic}{\operatorname{Pic}\nolimits}
\newcommand{\Spec}{\operatorname{Spec}\nolimits}
\newtheorem{theorem}{Theorem}[section]
\newtheorem{acknowledgement}[theorem]{Acknowledgement}
\newtheorem{algorithm}[theorem]{Algorithm}
\newtheorem{axiom}[theorem]{Axiom}
\newtheorem{case}[theorem]{Case}
\newtheorem{claim}[theorem]{Claim}
\newtheorem{conclusion}[theorem]{Conclusion}
\newtheorem{condition}[theorem]{Condition}
\newtheorem{conjecture}[theorem]{Conjecture}
\newtheorem{construction}[theorem]{Construction}
\newtheorem{corollary}[theorem]{Corollary}
\newtheorem{criterion}[theorem]{Criterion}
\newtheorem{definition}[theorem]{Definition}
\newtheorem{example}[theorem]{Example}
\newtheorem{exercise}[theorem]{Exercise}
\newtheorem{lemma}[theorem]{Lemma}
\newtheorem{notation}[theorem]{Notation}
\newtheorem{problem}[theorem]{Problem}
\newtheorem{proposition}[theorem]{Proposition}
\newtheorem{remark}[theorem]{Remark}
\newtheorem{solution}[theorem]{Solution}
\newtheorem{summary}[theorem]{Summary}
\newtheorem*{thm}{Theorem}

\renewcommand{\theequation}{\arabic{section}.\arabic{equation}}

\def \bp{{\mathbf p}}
\def \bA{{\mathbf A}}
\def \bL{{\mathbf L}}
\def \bF{{\mathbf F}}
\def \bS{{\mathbf S}}
\def \bC{{\mathbf C}}
\def \bD{{\mathbf D}}
\def \Z{{\Bbb Z}}
\def \F{{\Bbb F}}
\def \C{{\Bbb C}}
\def \N{{\Bbb N}}
\def \Q{{\Bbb Q}}
\def \G{{\Bbb G}}
\def \X{{\Bbb X}}
\def \P{{\Bbb P}}
\def \K{{\Bbb K}}
\def \E{{\Bbb E}}
\def \A{{\Bbb A}}
\def \BH{{\Bbb H}}
\def \T{{\Bbb T}}
\def\rH{{\mathrm{H}}}

\title[semi-derived Ringel-Hall algebras]{Semi-derived Ringel-Hall Algebras and Hall Algebras of Odd-Periodic Relative Derived Categories}
\thanks{Ji Lin is supported by the National Natural Science Foundation of China (Grant No. 12001107), the University Natural Science Project of Anhui Province(Grant No. KJ2021A0661), the University Outstanding Youth Research Project in Anhui Province (Grant No. 2022AH020082), the Scientific Research and Innovation Team Project of Fuyang Normal University (Grant No. TDJC2021009) and Liangang Peng is supported by the National Natural Science Foundation of China (Grant No. 11821001)}
\author[Lin]{Ji Lin}
\address{Department of Mathematics and Statistics, Fuyang Normal University, Fuyang 236037, P.R.China}
\email{jlin@fynu.edu.cn}
\author[Peng]{Liangang Peng}
\address{Department of Mathematics, Sichuan University, Chengdu 610064, P.R.China}
\email{penglg@scu.edu.cn}


\begin{abstract}
Let $t$ be a positive integer and $\ca$ a hereditary abelian category satisfying some finiteness conditions. We define the semi-derived Ringel-Hall algebra of $\ca$ from the category $\mathcal{C}_{\mathbb{Z}/t}(\mathcal{A})$ of $\mathbb{Z}/t$-graded complexes and obtain a natural basis of the semi-derived Ringel-Hall algebra. Moreover, we describe the semi-derived Ringel-Hall algebra by the generators and defining relations. In particular, if $t$ is an odd integer, we show that there is an embedding of derived Hall algebra of the odd-periodic relative derived category in the extended semi-derived Ringel-Hall algebra.
\end{abstract}
\maketitle
\tableofcontents
\section{Introduction}
Ringel-Hall algebras first introduced by Ringel \cite{R0} and generalized to the Kac-Moody type by Green \cite{Gr} have been not only a nice framework of the realization of quantum enveloping algebras, but also provide an important tool to encode the structure of categories (see \cite{Gr,R0,R2,R5} for more details).

Let $A$ be a finite-dimensional representation-finite hereditary algebra
over some finite field $k=\mathbb{F}_q$. The Ringel-Hall algebra of $A$ is an associative algebra over $\mathbb{Q}$ with a basis consisting of the isomorphism classes of finite-dimensional $A$-modules. A surprising result of Ringel stated that the positive part of quantum enveloping algebra of finite type can be obtained as the generic twisted Ringel-Hall algebra of $A$. Ringel also showed that the positive part of a simple Lie algebra can be recovered by the corresponding Ringel-Hall algebra.

After that much of the work mainly focused on how to recover the whole quantum enveloping algebras and the symmetrizable Kac-Moody Lie algebras in terms of Ringel-Hall algebras (see \cite{Br,Gr,P,P2,PX1,PX2,Rie,X}). Xiao \cite{X} successfully realized the whole quantum enveloping algebra via a direct model, i.e., the reduced Drinfeld double Ringel-Hall algebra by piecing together two Borel parts, but this is not intrinsic. Also L. Peng and Xiao constructed the theory of Ringel-Hall Lie algebras by which they recovered the corresponding simple Lie algebras as well as the symmetrizable Kac-Moody Lie algebras (see \cite{PX1, PX2} for more details).

On the other hand, To\"{e}n \cite{T} introduced the remarkable Hall algebra structure of a DG-enhanced triangulated category named derived Hall algebra. By using the octahedral axiom, Xiao and  Xu \cite{XX} directly showed To\"{e}n's formula for a triangulated category satisfying the left homological finiteness conditions. Xu and Chen \cite{XC} immediately generalized the work of To\"{e}n and Xiao and Xu to the odd-periodic triangulated category. Recently, Bridgeland \cite{Br} constructed the Hall algebra which is by definition the localization of the Ringel-Hall algebra of $\mathbb{Z}/2$-graded complexes with projective components at the set consisting of all quasi-isomorphisms. In particular, if the abelian category is hereditary, Bridgeland showed that the quantum enveloping algebra can be embedded into the so-called Bridgeland's Hall algebra. Moreover, Yanagida \cite{Ya16} showed that the Bridgeland's Hall algebra is isomorphic to the Drinfeld double of the twisted extended Ringel-Hall algebra. Zhang \cite{Zhang1} found the generators and relations of Bridgeland's Hall algebra of $\mathbb{Z}/t$-graded complexes over projectives for $t>2$ or $t=0$, here $\mathbb{Z}/0$-graded complex means the bounded complex. Zhang also realized the Lie algebra spanned by all indecomposable objects via Hall algebras of $\mathbb{Z}/1$-graded perfect complexes in \cite{Zhang2}.

As a generalization of Bridgeland's constuction, semi-derived Ringel-Hall algebras were introduced by Gorsky in terms of $\mathbb{Z}/2$-graded complexes over abelian categories with enough projective objects and Frobenius categories in \cite{Gor13} and \cite{Gor16} respectively. Afterwards, Sheng, Chen and Xu \cite{SCX} investigated the relations between the Hall algebra of an odd-periodic algebraic triangulated category and the semi-derived Hall algebra of the corresponding Frobenius category.

Inspired by the works of Bridgeland \cite{Br} and Gorsky \cite{Gor13, Gor16} on constructing Ringel-Hall algebras from the category of $\mathbb{Z}/2$-graded complexes, Lu and the second author \cite{LuP21} formulated semi-derived Ringel-Hall algebras over an arbitrary hereditary abelian category which may not have enough projective objects. It is shown that the semi-derived Ringel-Hall algebra is isomorphic to the Drinfeld double of the twisted extended Ringel-Hall algebra. Substituting the bounded complexes for the $\mathbb{Z}/2$-graded complexes, we explicitly gave an embedding of the derived Hall algebra in the semi-derived Ringel-Hall algebra by their generators in \cite{LinP}.

More recently, Lu and Wang \cite{LuW19,LuW20} introduced an $\imath$Hall algebra
which is by definition the twisted semi-derived Hall algebras for the $\imath$quiver algebras. They showed that $\imath$Hall algebras provide a realization of the quasi-split $\imath$quantum groups.

This paper is a sequel of our previous work in \cite{LuP21} and \cite{LinP}. In \cite{LinP} and the previous version of \cite{LuP21} (see \cite{LuP16}), the semi-derived Ringel-Hall algebras were called the {\it modified Ringel-Hall algebras}. To avoid confusion with the {\it modified quantum groups} (see \cite{Lu93}), we have changed the name and notation to the semi-derived Ringel-Hall algebra $\cs\cd\ch_{\mathbb{Z}/2}(\ca)$. In this paper, we mainly consider the semi-derived Ringel-Hall algebra of a hereditary abelian category over the category of $\mathbb{Z}/t$-graded complexes. However, it is not a trivial generalization of our previous work, especially for the case $t=1$. First of all, we modify the definition of semi-derived Ringel-Hall algebras further, which makes this kind of Hall algebras much more comprehensive than the one given in \cite{LuP21} for $t\neq1$.  As a result, we find an embedding of the derived Hall algebra in a ``larger'' semi-derived Ringel-Hall algebra which is by definition the tensor algebra of the semi-derived Ringel-Hall algebra and the extended quantum torus by adding ``the half elements'' to the quantum torus. And this construction is necessary for our embedding which is essentially different from the one given in \cite{SCX}.

The paper is organized as follows. In Section \ref{section locally homological finiteness of acyclic }, we first recall the relative derived categories following  \cite{Br, PX1}, and then investigate the locally homological finiteness of acyclic $\mathbb{Z}/t$-graded complexes. Consequently, we calculate the Euler forms for some special $\mathbb{Z}/t$-graded complexes. In Section \ref{Section definition and structure}, we define the semi-derived Ringel-Hall algebra in a slightly different way from the one we introduced in \cite{LuP21,LinP}. Next, we obtain a nice basis, as well as the generators and defining relations of the semi-derived Ringel-Hall algebras. In Section \ref{Hall algebras of odd-periodic relative derived categories}, we mainly recall the notion of derived Hall algebras of odd-periodic relative derived categories. Section \ref{main results} is devoted to establishing a relationship between the semi-derived Ringel-Hall algebras and the Hall algebras of odd-periodic relative derived categories. More specifically, we first form the extended semi-derived Ringel-Hall algebra and then directly construct an embedding of the derived Hall algebra in the extended semi-derived Ringel-Hall algebra via the generators.

\section{Locally homological finiteness of acyclic $\mathbb{Z}/t$-graded complexes}\label{section locally homological finiteness of acyclic }
In this section we mainly calculate Euler forms for some special periodic complexes. We first recall relative derived categories and refer to \cite{Br,PX1} for more details. Throughout this paper, $k$ denotes a finite field with $q$ elements and set $v:=\sqrt{q}$, and $|S|$ denotes the cardinality of a finite set $S$. Denote by $\ca$ an essentially small abelian $k$-linear category.
\subsection{Categories of $\mathbb{Z}/t$-graded complexes}
For any positive integer $t$, $\mathbb{Z}/t=\{0,1,\cdots, t-1\}$. A {\it $\mathbb{Z}/t$-graded complex} $M^\bullet=(M^i, d^i_{M^\bullet})_{i\in\mathbb{Z}/t}$ over $\ca$ consists of objects $M^i\in\ca$ and morphisms $d^i_{M^\bullet}: M^i\rightarrow M^{i+1}$ satisfying $d_{M^\bullet}^{i+1}d_{M^\bullet}^i=0(i\in\mathbb{Z}/t)$. Also it is convenient to express a $\mathbb{Z}/t$-graded complex by the following sequence
$$M^\bullet:= M^0\xrightarrow{d_{M^\bullet}^0}M^1\xrightarrow{d_{M^\bullet}^1}\cdots \xrightarrow{d_{M^\bullet}^{t-2}}M^{t-1}\xrightarrow{d_{M^\bullet}^{t-1}}M^0,$$
with $M^i\in\ca$ and $d_{M^\bullet}^{i+1}d_{M^\bullet}^i=0$ for any $i\in\mathbb{Z}/t$, where $M^i$ is the {\it $i$th homogenous component} of $M^\bullet$ and $d_{M^\bullet}^i$ is the {\it $i$th differential} of $M^\bullet$. The \emph{$i$th relative homology group} $H^i(M^\bullet)$ of $M^\bullet$ is defined to be $$H^i(M^\bullet)=\Ker(d_{M^\bullet}^i)/\Im(d_{M^\bullet}^{i-1})$$
for any $i\in\mathbb{Z}/t$. And it is {\it acyclic} if $H^i(M^\bullet)=0$ for each $i\in\mathbb{Z}/t$.

Let $\cc_{\mathbb{Z}/t}(\ca)$ be the {\it category of $\mathbb{Z}/t$-graded complexes} over $\ca$. A morphism $$f^\bullet=(f^i)_{i\in\mathbb{Z}/t}: M^\bullet=(M^i, d^i_{M^\bullet})_{i\in\mathbb{Z}/t}\rightarrow N^\bullet=(N^i, d_{N^\bullet}^i)_{i\in\mathbb{Z}/t}$$ of $\mathbb{Z}/t$-graded complexes is a sequence of morphisms $f^i: M^i\rightarrow N^i$ satisfying $d_{N^\bullet}^{i}f^i=f^{i+1}d_{M^\bullet}^i$ for any $i\in\mathbb{Z}/t$. Two morphisms $f^\bullet, g^\bullet: M^\bullet\rightarrow N^\bullet$ are {\it relatively homotopic} if there are morphisms $s^i: M^i\rightarrow N^{i-1}$ such that $$f^i-g^i=s^{i+1}d^i_{M^\bullet}+d^{i-1}_{N^\bullet}s^{i}$$ for any $i\in\mathbb{Z}/t$. This is an equivalence relation and denoted by $\ck_{\mathbb{Z}/t}(\ca)$ the induced {\it relative homotopy category}. In addition, there is a shift functor of $\cc_{\mathbb{Z}/t}(\ca)$(resp. $\ck_{\mathbb{Z}/t}(\ca)$) denoted by $[1]$ which defines a new complex $M^\bullet[1]$ for any $M^\bullet=(M^i, d_{M^\bullet}^i)\in\cc_{\mathbb{Z}/t}(\ca)$ as follows:
$$(M^\bullet[1])^{i}=M^{i+1}~\mbox{and}~d_{M^\bullet[1]}^i=-d_{M^\bullet}^{i+1}$$
for any $i\in\mathbb{Z}/t$ and it induces an automorphism of $\cc_{\mathbb{Z}/t}(\ca)$(resp. $\ck_{\mathbb{Z}/t}(\ca)$) in an obvious way. A morphism $f^\bullet: M^\bullet\rightarrow N^\bullet$ in $\cc_{\mathbb{Z}/t}(\ca)$ is a {\it quasi-isomorphism} if the induced morphisms $$H^i(f^\bullet): H^i(M^\bullet)\longrightarrow H^i(N^\bullet)$$ are isomorphisms for all $i\in\mathbb{Z}/t$.
The {\it relative derived category} $\cd_{\mathbb{Z}/t}(\ca)$ is the localization of $\ck_{\mathbb{Z}/t}(\ca)$ with respect to the set of all the quasi-isomorphisms. Then one can get that both $\ck_{\mathbb{Z}/t}(\ca)$ and $\cd_{\mathbb{Z}/t}(\ca)$ are triangulated categories with the translation functor $T=[1]$.

In particular, $\mathbb{Z}/1$-graded complex denoted by $M^\bullet=(M, d_M)$ consists of an object $M\in\ca$ and an endomorphism $d_M: M\rightarrow M$ satisfying $d_M^2=0$. Moreover, one can easily obtain an exact functor
$$C_1: \cc_{\mathbb{Z}/t}(\ca)\rightarrow \cc_{\mathbb{Z}/1}(\ca)$$
defined by $$M^\bullet\mapsto C_1(M^{\bullet}):=\left(\bigoplus_{i=0}^{t-1}M^{i}, \left(\begin{array}{ccccc}
0&0&\cdots&0&d_{M^\bullet}^{t-1}\\
d_{M^\bullet}^0&0&\cdots&0&0\\
0&d_{M^\bullet}^1&\cdots&0&0\\
\multicolumn{5}{c}{\dotfill}\\
0&0&\cdots&d_{M^\bullet}^{t-2}&0
\end{array}\right)\right)$$
for any $M^\bullet:= M^0\xrightarrow{d_{M^\bullet}^0}M^1\xrightarrow{d_{M^\bullet}^1}\cdots \xrightarrow{d_{M^\bullet}^{t-2}}M^{t-1}\xrightarrow{d_{M^\bullet}^{t-1}}M^0$ in $\cc_{\mathbb{Z}/t}(\ca)$ in case $t>1$. Conversely, one can also form a $\mathbb{Z}/t$-graded complex
$$C_t(M^{\bullet}):=M\xrightarrow{d_{M}}M\xrightarrow{d_{M}}\cdots\xrightarrow{d_{M}} M\xrightarrow{d_{M}}M$$
from any $\mathbb{Z}/1$-graded complex $M^{\bullet}=(M, d_{M})\in\cc_{\mathbb{Z}/1}(\ca)$. Obviously, $C_t(M^{\bullet})$ is acyclic if and only if $M^{\bullet}$ is acyclic.

\begin{lemma}\label{lemma of split complex}
Let $M^{\bullet}=(M, d_M)$ be a complex in $\cc_{\mathbb{Z}/1}(\ca)$. Then
\begin{itemize}
\item[(1)] $M^{\bullet}$ is a direct summand of $C_1\left(C_t(M^{\bullet})\right)$ for some $t\geq 2$;\\
\item[(2)] $C_t(M^{\bullet})\in\cc_{\mathbb{Z}/t}(\ca)$ and $C_1(C_t(M^{\bullet}))\in\cc_{\mathbb{Z}/1}(\ca)$ are both acyclic complexes if $M^{\bullet}$ is acyclic.
\end{itemize}
\end{lemma}
\begin{proof}
(1) Assume that $k$ is of characteristic $p$. Let $t=p+1$, and then the following morphism
$$M^\bullet\xrightarrow{\left(\begin{array}{cccc}1&1&\cdots&1\end{array}\right)^{T}}C_1(C_t(M^\bullet))$$
is a section since it admits a left inverse $$C_1(C_t(M^\bullet))\xrightarrow{\left(\begin{array}{cccc}\frac{1}{t}&\frac{1}{t}&\cdots&\frac{1}{t}\end{array}\right)}M^\bullet.$$
(2) This is straightforward.
\end{proof}

\subsection{Locally homological finiteness of acyclic $\mathbb{Z}/t$-graded complexes}
Given an object $X\in\ca$ and $i\in\mathbb{Z}/t$, denote by $K_{X, i}$ the following acyclic $\mathbb{Z}/t$-graded complex
$$0\rightarrow\cdots\rightarrow 0\rightarrow X\xrightarrow{1_X}X\rightarrow 0\rightarrow\cdots \rightarrow 0$$
with $X$ sitting in degrees $i-1$ and $i$ in case $t>1$, and let $K_{X}:=K_{X,0}$ be the acyclic $\mathbb{Z}/1$-graded complex $\left(X\oplus X,\begin{pmatrix}
0&1\\0&0
\end{pmatrix}\right)$. The following $\mathbb{Z}/t$-graded stalk complex
$$0\rightarrow\cdots\rightarrow 0\rightarrow X \rightarrow 0\rightarrow\cdots \rightarrow 0$$
with $X$ concentrated in degree $i$ is denoted by $U_{X, i}$, also set $U_X:=U_{X,0}$ in case $t=1$.

Assume now that $\ca$ is, in addition, hereditary. And then we have the following result.

\begin{proposition}\label{proposition t-periodic iso to the homology}
Let $M^{\bullet}=(M^i, d^i_{M^{\bullet}})_{i\in\mathbb{Z}/t}$ be a $\mathbb{Z}/t$-graded complex in $\cc_{\mathbb{Z}/t}(\ca)$. Then $M^\bullet$ is isomorphic to $\bigoplus\limits_{i=0}^{t-1}U_{H^i(M^\bullet),i}$ in $\cd_{\mathbb{Z}/t}(\ca)$.
\end{proposition}
\begin{proof}
Since $\ca$ is hereditary we know that the action of $\Ext^1_{\ca}(H^i(M^\bullet),-)$ on the canonical
epimorphism $M^{i-1}\twoheadrightarrow \Im(d^{i-1}_{M^\bullet})$ induces an epimorphism $\Ext_{\ca}^1(H^i(M^\bullet),M^{i-1})\twoheadrightarrow
\Ext_{\ca}^1(H^i(M^\bullet),\Im(d^{i-1}_{M^\bullet}))$, there exists object $B^i(i\in\mathbb{Z}/t)$,  morphisms $l^i: M^i\rightarrow B^{i+1}$ and $g^i: B^i\rightarrow\Ker(d^i_{M^\bullet})$ satisfying the following commutative diagram
\[\xymatrix{0\ar[r]&M^{i-1}\ar[r]^{l^{i-1}}\ar[d]^{\overline{d^{i-1}_{M^{\bullet}}}}&B^i\ar[r]^{p^i g^i}\ar[d]^{g^i}&H^i(M^\bullet)\ar@{=}[d]\ar[r]&0\\
0\ar[r]&\Im(d_{M^{\bullet}}^{i-1})\ar[r]^{j_1^i}&
\Ker(d_{M^{\bullet}}^{i})\ar[r]^{p^i}&H^i(M^\bullet)\ar[r]&0.}\]
So we have the following commutative diagram
\[
\xymatrix{&0\ar[d]&0\ar[d]\\
&\bigoplus\limits_{i=0}^{t-1}U_{\Ker(d_{M^{\bullet}}^{i-1}),i}\ar@{=}[r]\ar[d]^{j_2^\bullet}&\bigoplus\limits_{i=0}^{t-1}U_{\Ker(d_{M^{\bullet}}^{i-1}),i}\ar[d]^{l^\bullet j_2^\bullet}&\\
0\ar[r]&M^\bullet[-1]\ar[r]^{l^\bullet}\ar[d]^{\overline{d^\bullet_{M^{\bullet}}}}&B^\bullet\ar[r]^{p^\bullet g^\bullet}\ar[d]^{g^\bullet}&\bigoplus\limits_{i=0}^{t-1}U_{H^i(M^\bullet),i}\ar@{=}[d]\ar[r]&0\\
0\ar[r]&\bigoplus\limits_{i=0}^{t-1}U_{\Im(d_{M^{\bullet}}^{i-1}),i}\ar[d]\ar[r]^{j_1^\bullet}&\bigoplus\limits_{i=0}^{t-1}U_{\Ker(d_{M^{\bullet}}^{i}),i}\ar[d]\ar[r]^{p^\bullet}&\bigoplus\limits_{i=0}^{t-1}U_{H^i(M^\bullet),i}\ar[r]&0\\
&0&0}
\]
in the category $gr_{\mathbb{Z}/t}(\ca)$ of graded objects, whose objects are $\mathbb{Z}/t$-graded families of objects of $\ca$ and morphisms are given componentwise, where all the rows and columns are short exact sequences, and $l^\bullet=diag(l^{t-1},l^0,\cdots,l^{t-2})$, $\overline{d^\bullet_{M^{\bullet}}}=diag(\overline{d_{M^{\bullet}}^{t-1}},\overline{d_{M^{\bullet}}^{0}},\cdots,\overline{d_{M^{\bullet}}^{t-2}}),j_1^\bullet=diag(j_1^{t-1},j_1^{0},\cdots,j_1^{t-2})$ and $j_2^\bullet=diag(j_2^{t-1},j_2^{0},\cdots,j_2^{t-2})$ are the homomorphisms naturally induced by the differentials of $M^\bullet$ satisfying $d_{M^{\bullet}}^i=j_2^{i+1}j_1^i\overline{d_{M^{\bullet}}^i}$. Let $d^i_{B^\bullet}=-l^{i} j_2^{i}g^i(i\in \Z/t)$ be the differential of $B^\bullet$, then one can easily check that $B^\bullet$ is acyclic.

Furthermore, one can obtain the following short exact sequence $$0\rightarrow M^\bullet[-1]\xrightarrow{l^\bullet} B^\bullet\xrightarrow{p^\bullet g^\bullet }\bigoplus\limits_{i=0}^{t-1}U_{H^i(M^\bullet),i}\rightarrow 0$$
in $\cc_{\Z/t}(\ca)$ which implies that $M^\bullet$ is isomorphic to $\bigoplus\limits_{i=0}^{t-1}U_{H^i(M^\bullet),i}$ in $\cd_{\Z/t}(\ca)$ as required.
\end{proof}

The following important lemma
plays a key role in the locally homological finiteness of acyclic $\mathbb{Z}/t$-graded complexes.

\begin{lemma}\label{lemma homological finiteness K_{X,i}}
Let $X$ be an object in $\ca$. Then the Ext-dimension of $K_{X,i}(i\in\mathbb{Z}/t)$ is less than or equal to 1, i.e., $$\Ext^p_{\cc_{\Z/t}(\ca)}(K_{X, i}, M^\bullet)=0\ \text{and}\ \Ext^p_{\cc_{\Z/t}(\ca)}(M^\bullet, K_{X, i})=0$$ for any $M^\bullet\in\cc_{\Z/t}(\ca)$ and $p\geq 2$.
\end{lemma}

\begin{proof}
It suffices to show the lemma in case $t=1$, since a similar proof in \cite{LuP21} yields the above identities if $t>1$.

For any short exact sequence $\xi^\bullet$ in $\Ext^2_{\cc_{\Z/1}(\ca)}(K_{X}, M^\bullet)$, assume that it is of the following form
$$0\rightarrow M^\bullet \xrightarrow{f^\bullet_{3}} M_2^\bullet\xrightarrow{f^\bullet_{2}} M_1^\bullet\xrightarrow{f^\bullet_{1}}K_{X}\rightarrow 0,$$
where $M^\bullet=(M,d)$, $M_2^\bullet=(M_2,d_2)$, $M_1^\bullet=(M_1,d_1)$. Note that $K_X=\left(X\oplus X,\begin{pmatrix}
0&1\\0&0
\end{pmatrix}\right)$ yields that $f^\bullet_{1}=\begin{pmatrix}f_{11}\\f_{11}d_1\end{pmatrix}$ for some epimorphism $f_{11}\in \Hom_{\ca}(M_1,X)$. So we can form the pull-back of $\begin{pmatrix}f_{11}\\f_{11}d_1\end{pmatrix}$ and $\begin{pmatrix}0\\1\end{pmatrix}$ in $\ca$, hence there exists an object $W_1\in\ca$, an epimorphism $g_1\in\Hom_{\ca}(W_1,X)$ and a homomorphism $h_1\in\Hom_{\ca}(W_1,M_1)$ satisfying the following commutative diagram
\[
\xymatrix{
W_1\ar@{->>}[r]^{g_1}\ar[d]_{h_1}&X\ar[d]^{\begin{pmatrix}0\\1\end{pmatrix}}\\
M_1\ar@{->>}[r]_{\begin{pmatrix}f_{11}\\f_{11}d_1\end{pmatrix}}&X\oplus X,
}
\]
which implies the next commutative diagram
\[
\xymatrix{
0\ar[r]&M\ar[r]^{f_3}\ar@{=}[d]&M_2\ar[r]^{g_2}\ar@{=}[d]&W_1\ar[r]^{g_1}\ar[d]^{h_1}&X\ar[r]\ar[d]^{\begin{pmatrix}0\\1\end{pmatrix}}&0\\
0\ar[r]&M^\bullet\ar[r]^{f_3^\bullet}&M_2^\bullet\ar[r]^{f_2^\bullet}&M_1^\bullet\ar[r]^{f_1^\bullet}&K_X\ar[r]&0
}\]
in $\ca$ with exact rows. Furthermore, it is not hard to obtain the following exact sequence
$$\eta^\bullet: 0\rightarrow K_{M}\xrightarrow{g_3^\bullet}K_{M_2}\xrightarrow{g_2^\bullet}K_{W_1}\xrightarrow{g_1^\bullet}K_{X}\rightarrow 0$$in $\cc_{\Z/1}(\ca)$, where $g_3^\bullet=\begin{pmatrix}f_3&0\\0&f_3\end{pmatrix}$, $g_2^\bullet=\begin{pmatrix}g_2&0\\0&g_2\end{pmatrix}$,  $g_1^\bullet=\begin{pmatrix}g_1&0\\0&g_1\end{pmatrix}$. In addition one can also get a morphism from $\eta^\bullet$ to $\xi^\bullet$ as follows
\[
\xymatrix{
0\ar[r]&K_M\ar[r]^{g_3^\bullet}\ar[d]^{h_3^\bullet}&K_{M_2}\ar[r]^{g_2^\bullet}\ar[d]^{h_2^\bullet}&K_{W_1}\ar[r]^{g_1^\bullet}\ar[d]^{h^\bullet}&K_X\ar[r]\ar@{=}[d]&0\\
0\ar[r]&M^\bullet\ar[r]^{f_3^\bullet}&M_2^\bullet\ar[r]^{f_2^\bullet}&M_1^\bullet\ar[r]^{f_1^\bullet}&K_X\ar[r]&0,
}\]
which means that $\Ext^2_{\cc_{\Z/1}(\ca)}(K_X, h_3^\bullet)([\eta^\bullet])=[\xi^\bullet]$, where $h_3^\bullet=(d,1)$, $h_2^\bullet=(d_2,1)$, $h^\bullet=\begin{pmatrix}d_1h_1,h_1\end{pmatrix}$. We get that $[\eta^\bullet]=0$ since $\ca$ is hereditary, and then $[\xi^\bullet]=0$. Therefore, $\Ext^p_{\cc_{\Z/1}(\ca)}(K_{X}, M^\bullet)=0$ as required is directly induced by Yoneda product. The proof of the second equality is by duality.
\end{proof}

\begin{proposition}\label{proposition homo-fin acyclic t-periodic complexes}
Let $M^\bullet, N^\bullet$ be two complexes in $\cc_{\Z/t}(\ca)$. If $N^\bullet$ is acyclic, then
$$\Ext^p_{\cc_{\Z/t}(\ca)}(N^\bullet, M^\bullet)=0\ \text{and}\ \Ext^p_{\cc_{\Z/t}(\ca)}(M^\bullet, N^\bullet)=0$$ for any $M^\bullet\in\cc_{\Z/t}(\ca)$ and  $p\geq 2$.
\end{proposition}
\begin{proof}
As shown in \cite{LP}, it is easy to get that the Galois covering functor $F_t: \cd^b(\ca)\rightarrow \cd_{\Z/t}(\ca)$ is dense in case $t>1$, so it immediately follows from a similar proof in \cite{LuP21} in case $t>1$. And then we just need to show the above two identities in case $t=1$.

For any acyclic complex $N^\bullet=(N,d)$ in $\cc_{\Z/1}(\ca)$, we only prove the first identity and the proof of the other one is by duality. By the similar proof of Lemma 2.1 in \cite{LuP21}, there exists a bounded complex $V^\bullet$ satisfying the following short exact sequence of $\mathbb{Z}/t$-graded complexes
$$0\rightarrow C_t(N^\bullet)\xrightarrow{qis} F_t(V^\bullet)\rightarrow K_{N/\Im(d),1}\rightarrow 0,$$
which yields the next short exact sequence
$$0\rightarrow C_1(C_t(N^\bullet))\xrightarrow{qis}C_1(F_t(V^\bullet))\rightarrow K_{N/\Im(d)}\rightarrow 0.$$
Moreover, $F_t(V^\bullet)$ can be obtained from the acyclic complexes of the form $K_{X,i}$ for some $X\in\ca$ and $i\in \mathbb{Z}/t$ via iterated extensions. Hence the induced acyclic $\mathbb{Z}/1$-graded complex $C_1(F_t(V^\bullet))$ can be also obtained from the acyclic complexes of the form $K_{X}$ for some $X\in\ca$ by iterated extensions. It means that $\Ext^p_{\cc_{\Z/1}(\ca)}(C_1(F_t(V^\bullet)),M^\bullet)=0$ for any $M^\bullet\in\cc_{\Z/1}(\ca)$ and $p\geq 2$ by Lemma \ref{lemma homological finiteness K_{X,i}}. So one can easily deduce that $\Ext^p_{\cc_{\Z/1}(\ca)}(C_1(C_t(N^\bullet)), M^\bullet)=0$. And it follows from Lemma \ref{lemma of split complex} that $N^\bullet$ is a direct summand of $C_1(C_t(N^\bullet))$ for some $t>1$, hence $\Ext^p_{\cc_{\Z/1}(\ca)}(N^\bullet, M^\bullet)=0$.
This finishes the proof.
\end{proof}
The above result is necessary for the realization of the universal $q$-Onsager algebra via $\imath$Hall algebra of the projective line over a finite field, and there is a proof in \cite{LRW20} following our method.
\subsection{The Euler forms}\label{section Euler form}
Assume now that $\ca$ satisfies the following two additional finiteness conditions $$\dim_k\Hom_{\ca}(X, Y)<\infty~\text{and}~\dim_{k}\Ext^1_{\ca}(X, Y)<\infty$$ for any objects $X, Y$ in $\ca$. Let $\Iso(\ca)$ denote the isomorphism classes of $\ca$, and for any $X\in\ca$ the image in the Grothendieck group $K_0(\ca)$ is denoted by $\widehat{X}$. Since $\ca$ is hereditary, the multiplicative {\it Euler form}
$$\langle-, -\rangle: K_0(\ca)\times K_0(\ca)\rightarrow \mathbb{Q}^\times$$
defined by $$\langle\widehat{X}, \widehat{Y}\rangle=\frac{|\Hom_{\ca}(X, Y)|}{|\Ext^1_{\ca}(X, Y)|}$$
for $X, Y\in\ca$ is a bilinear form on $K_0(\ca)$, and the {\it symmetric Euler form} is defined to be $$(\widehat{X}, \widehat{Y}):=\langle\widehat{X}, \widehat{Y}\rangle\langle\widehat{Y}, \widehat{X}\rangle.$$

Let $\cc_{\Z/t}^{ac}(\ca)$ be the subcategory of $\cc_{\Z/t}(\ca)$ formed by all the acyclic complexes in $\cc_{\Z/t}(\ca)$. By Proposition \ref{proposition homo-fin acyclic t-periodic complexes} we can define the {\em Euler forms}
$$\langle-, -\rangle: \Iso(\cc_{\Z/t}^{ac}(\ca))\times \Iso(\cc_{\Z/t}(\ca))\rightarrow \mathbb{Q}^\times$$ and $$\langle-, -\rangle: \Iso(\cc_{\Z/t}(\ca))\times \Iso(\cc_{\Z/t}^{ac}(\ca))\rightarrow \mathbb{Q}^\times$$
by
$$\langle[K^{\bullet}], [M^{\bullet}]\rangle=\prod_{p=0}^{+\infty}|\Ext^p_{\cc_{\Z/t}(\ca)}(K^{\bullet}, M^{\bullet})|^{(-1)^p}=\frac{|\Hom_{\cc_{\Z/t}(\ca)}(K^{\bullet}, M^{\bullet})|}{|\Ext^1_{\cc_{\Z/t}(\ca)}(K^{\bullet}, M^{\bullet})|}$$
and
$$\langle[M^{\bullet}], [K^{\bullet}]\rangle=\prod_{p=0}^{+\infty}|\Ext^p_{\cc_{\Z/t}(\ca)}(M^{\bullet}, K^{\bullet})|^{(-1)^p}=\frac{|\Hom_{\cc_{\Z/t}(\ca)}(M^{\bullet}, K^{\bullet})|}{|\Ext^1_{\cc_{\Z/t}(\ca)}(M^{\bullet}, K^{\bullet})|}$$
for $[K^{\bullet}], [M^{\bullet}]\in\Iso(\cc_{\Z/t}(\ca))$ with $K^{\bullet}$ acyclic. Obviously, they can both descend to the bilinear forms on the Grothendieck groups $K_0(\cc_{\Z/t}(\ca))$ and $K_0(\cc_{\Z/t}^{ac}(\ca))$, also denoted by the same symbols:
$$\langle-,-\rangle: K_0(\cc_{\Z/t}^{ac}(\ca))\times K_0(\cc_{\Z/t}(\ca))\longrightarrow \mathbb{Q^{\times}}$$
and
$$\langle-,-\rangle: K_0(\cc_{\Z/t}(\ca))\times K_0(\cc_{\Z/t}^{ac}(\ca))\longrightarrow \mathbb{Q^{\times}}.$$

\begin{proposition}\label{prop euler form of acycli complexes}
For any $X, Y\in\ca$, we have

\begin{itemize}
\item [(1)] $\langle[K_{X, i}], [U_{Y,j}]\rangle=\begin{cases}\langle\widehat{X}, \widehat{Y}\rangle^{\delta_{t-1}^j}& i=0,\\ \langle\widehat{X}, \widehat{Y}\rangle^{\delta_{i-1}^j}& i=1,\cdots,t-1,
    \end{cases}$
     and $\langle[U_{X,i}], [K_{Y, j}]\rangle=\langle\widehat{X}, \widehat{Y}\rangle^{\delta^{j}_i}$,
\item [(2)] $\langle[K_{X, i}][K_{Y,j}]\rangle=\begin{cases} \langle\widehat{X}, \widehat{Y}\rangle^{(\delta_{0}^j+\delta_{t-1}^j)} & i=0, \\ \langle\widehat{X}, \widehat{Y}\rangle^{(\delta_{i}^j+\delta_{i-1}^j)}& i=1,\cdots,t-1,\end{cases}$\\
\end{itemize}for each $i,j\in\Z/t$, in case $t>1$. And if $t=1$, then
\begin{itemize}
\item[(3)] $\langle [K_X], [U_Y]\rangle=\langle\widehat{X}, \widehat{Y}\rangle$ and $\langle[U_X], [K_Y]\rangle=\langle\widehat{X}, \widehat{Y}\rangle$,
\item[(4)] $\langle [K_X], [K_Y]\rangle=\langle\widehat{X}, \widehat{Y}\rangle^2$.
\end{itemize}
\end{proposition}
\begin{proof}
$(1)$ and (2) can be directly deduced from the similar calculations as in \cite{LuP21, LinP}.

$(3)$ We just show the first identity and omit the similar proof of the other one. Since one can easily get the isomorphism $$\Hom_{\cc_{\Z/1}(\ca)}(K_X, U_Y)\cong \Hom_{\ca}(X, Y)\cong\Hom_{\cc_{\Z/1}(\ca)}(U_X,U_Y),$$ it suffices to show that $$|\Ext^1_{\cc_{\Z/1}(\ca)}(K_X, U_Y)|=|\Ext^1_{\ca}(X, Y)|.$$

Note that we have the following exact sequence\\
$0\rightarrow \Hom_{\cc_{\Z/1}(\ca)}(U_X,U_Y)\rightarrow\Ext^1_{\cc_{\Z/1}(\ca)}(U_X,U_Y)\xrightarrow{\psi} \Ext^1_{\cc_{\Z/1}(\ca)}(K_X,U_Y)\rightarrow \Ext^1_{\cc_{\Z/1}(\ca)}(U_X,U_Y)\rightarrow\Ext^2_{\cc_{\Z/1}(\ca)}(U_X,U_Y)\rightarrow 0$\\
induced by the short exact sequence $0\rightarrow U_X\xrightarrow{\begin{pmatrix}1\\0\end{pmatrix}}K_{X}\xrightarrow{(0,1)}U_X\rightarrow0$. And we claim that $\psi=\Ext^1_{\cc_{\Z/1}(\ca)}((0,1),U_Y)$ is an epimorphism. To this end, for any element $[\xi^\bullet]$ of $\Ext^1_{\cc_{\Z/1}(\ca)}(K_X,U_Y)$, without loss of generality, assume that $\xi^\bullet$ is a short exact sequence of the following form
$$\xi^\bullet: 0\rightarrow U_Y\xrightarrow{f}M^\bullet\xrightarrow{\left(\begin{array}{c}g_1\\g_2\end{array}\right)}K_{X}\rightarrow 0,$$
for some $M^\bullet=(M, d_{M})\in\cc_{\Z/1}(\ca)$ . Then it is not hard to see that there exists a homomorphism $\left(h_1, h_2\right): K_X\rightarrow M^\bullet$ in $\cc_{\Z/1}(\ca)$ such that the following diagram is  commutative
\[
\xymatrix{
0\arrow[r]&U_Y\arrow[r]^{f}\arrow[d]_{0}&M^\bullet\arrow[r]^{\left(\begin{array}{c}g_1\\g_2\end{array}\right)}\arrow[d]_{d_M}&~~K_X~~\arrow[r]\ar@{-->}[ld]_{\tiny \left(h_1, h_2\right)}\arrow[d]^{\tiny \left(\begin{array}{cc}0&1\\0&0\end{array}\right)}&0\\
0\arrow[r]&U_Y\arrow[r]_{f}&M^\bullet\arrow[r]_{\left(\begin{array}{c}g_1\\g_2\end{array}\right)}&~~K_X~~\arrow[r]&0,}
\]
hence $d_M=\left(h_1, h_2\right)\left(\begin{array}{c}g_1\\g_2\end{array}\right)$ and $\left(\begin{array}{cc}0&1\\0&0\end{array}\right)=\left(\begin{array}{cc}g_1h_1&g_1h_2\\g_2h_1&g_2h_2\end{array}\right)$.
As a result, one can get that $h_1g_2=0$ since $d_M^2=0$. Thus $h_1=0$ and $d_M=h_2g_2$. Consequently, one can also obtain the following cummutative diagram in $\cc_{\Z/1}(\ca)$

\[\xymatrix{
&&0\ar[d]&0\ar[d]&\\
&&U_Y\ar@{=}[r]\ar[d]^{f}&U_Y\ar@{-->}[d]^{i}&\\
0\ar[r]&U_X\ar@{=}[d]\ar[r]^{h_2}&M^\bullet\ar[d]^{\tiny\left(\begin{array}{c}g_1\\g_2\end{array}\right)}\ar[r]^{q}&L^\bullet\ar[r]\ar@{-->}[d]^{p}&0\\
0\ar[r]&U_X\ar[r]_{\tiny \left(\begin{array}{c}1\\0\end{array}\right)}&K_X\ar[d]\ar[r]_{(0, 1)}&U_X\ar[d]\ar[r]&0.\\
&&0&0
}\]

Thus we have the following exact sequence
$$\eta^\bullet: 0\rightarrow U_Y\xrightarrow{i}L^\bullet\xrightarrow{p}U_X\rightarrow 0,$$
and then it is easy to prove that $\psi([\eta^\bullet])=[\xi^\bullet]$ by using the above commutative diagram. So we get the following short exact sequence
$$0\rightarrow \Hom_{\cc_{\Z/1}(\ca)}(U_X,U_Y)\rightarrow\Ext^1_{\cc_{\Z/1}(\ca)}(U_X,U_Y)\xrightarrow{\psi} \Ext^1_{\cc_{\Z/1}(\ca)}(K_X,U_Y)\rightarrow 0.$$

On the other hand, one can easily obtain the following epimorphism
$$\theta:\Ext^1_{\cc_{\Z/1}(\ca)}(U_X,U_Y)\twoheadrightarrow \Ext^1_{\ca}(X, Y)$$
with $\Ker(\theta)\cong\Hom_{\ca}(X, Y)$. So  $|\Ext^1_{\cc_{\Z/1}(\ca)}(K_X, U_Y)|=|\Ext^1_{\ca}(X, Y)|.$

$(4)$ follows immediately from the identity $\widehat{K_X}=\widehat{U_X}+\widehat{U_X}$ in $K_0(\cc_{\Z/1}(\ca))$ and $(3)$.
\end{proof}

\begin{proposition}\label{proposition euler form of acyclic and any complex}
Let $K^\bullet=(K^i, d^i_{K^\bullet})_{i\in\Z/t}\in\cc_{\Z/t}(\ca)$ be an acyclic complex. Then we have
$$\langle [K^\bullet], [M^\bullet]\rangle=\prod_{i\in\Z/t}\langle [K_{\Im(d_{K^\bullet}^i),i+1}], [M^\bullet]\rangle\ \ \text{and}\ \ \langle [M^\bullet], [K^\bullet]\rangle=\prod_{i\in\Z/t}\langle [M^\bullet], [K_{\Im(d_{K^\bullet}^i),i+1}]\rangle$$
for any complex $M^\bullet=(M^i, d_{M^\bullet}^i)_{i\in\Z/t}\in\cc_{\Z/t}(\ca)$.
\end{proposition}
\begin{proof}
In case $t>1$, there is a bounded complex
$$U^\bullet: \cdots\rightarrow 0\rightarrow K^0\xrightarrow{g^0}L\xrightarrow{d_{K^\bullet}^1h^1}K^2\xrightarrow{d_{K^\bullet}^2}\cdots K^{t-1}\xrightarrow{d^{t-1}}\Ker(d_{K^\bullet}^0)\rightarrow 0\rightarrow\cdots$$
induced by the following commutative diagram with exact rows and columns
\[\xymatrix{&0\ar[d]&0\ar[d]\\
&\Ker(d_{K^\bullet}^0)\ar@{=}[r]\ar[d]^{i^0}&\Ker(d_{K^\bullet}^0)\ar[d]^{i^1}&\\
0\ar[r]&K^0\ar[r]^{g^0}\ar[d]^{h^0}&L\ar[r]\ar[d]^{h^1}&\Coker(d_{K^\bullet}^0)\ar@{=}[d]\ar[r]&0\\
0\ar[r]&\Im(d_{K^\bullet}^0)\ar[d]\ar[r]^{g^1}&K^1\ar[d]\ar[r]&\Coker(d_{K^\bullet}^0)\ar[r]&0.\\
&0&0}
\]
And then we obtain the following short exact sequence in $\cc_{\Z/t}(\ca)$
$$0\rightarrow K_{\Ker(d_{K^\bullet}^0),1}\rightarrow F_t(U^\bullet)\rightarrow K^\bullet\rightarrow 0,$$
where $F_t: \cd^b(\ca)\rightarrow \cd_{\Z/t}(\ca)$ denotes the Galois covering functor and $$F_t(U^\bullet):= K^0\oplus\Ker(d_{K^\bullet}^0)\xrightarrow{\left(g^0, 0\right)}L\xrightarrow{d_{K^\bullet}^1h^1}K^2\rightarrow\cdots\rightarrow K^{t-1}\xrightarrow{\left(\begin{array}{c}0\\d_{K^\bullet}^{t-1}\end{array}\right)}K^0\oplus\Ker(d_{K^\bullet}^0).$$
We can also form the following exact sequence in $\cc_{\Z/t}(\ca)$
$$0\rightarrow K_{K^0,1}\rightarrow F_t(U^\bullet)\rightarrow N^\bullet\rightarrow 0,$$
where $N^\bullet=\Ker(d_{K^\bullet}^0)\xrightarrow{0}\Coker(d_{K^\bullet}^0)\xrightarrow{d^1}K^2
\rightarrow\cdots\rightarrow K^{t-1}\xrightarrow{d_{K^\bullet}^{t-1}}\Ker(d_{K^\bullet}^0)$. By \cite[Lemma 3.5]{LinP} we know that
$$\widehat{N^\bullet}=\widehat{K_{\Im(d_{K^\bullet}^1),2}}+\widehat{K_{\Im(d_{K^\bullet}^2),3}}+\cdots+\widehat{K_{\Im(d_{K^\bullet}^{t-1}),0}},$$ thus one can get the required identity $\langle [K^\bullet], [M^\bullet]\rangle=\prod_{i\in\Z/t}\langle [K_{\Im(d_{K^\bullet}^i),i+1}], [M^\bullet]\rangle$.

In case $t=1$, we denote by $K^\bullet=(K, d_{K})$ and $M^\bullet=(M, d_{M})$ the corresponding $\mathbb{Z}/1$-graded complexes. By definition one can easily get that
\begin{eqnarray*}
\langle [K^\bullet], [U_{\Ker(d_M)}]\rangle^{2} &=&\langle[K^\bullet], [K_{\Ker(d_M)}]\rangle\\
&=&\langle[U_{\Im(d_K)}], [K_{\Ker(d_M)}]\rangle^{2}\\
&=&\langle[K_{\Im(d_K)}], [K_{\Ker(d_M)}]\rangle\\
&=&\langle[K_{\Im(d_K)}], [U_{\Ker(d_M)}]\rangle^{2}.
\end{eqnarray*}
Similarly, we can also obtain that $\langle [K^\bullet], [U_{\Im(d_M)}]\rangle^{2}=\langle[K_{\Im(d_K)}], [U_{\Im(d_M)}]\rangle^{2}$.
Therefore, we get that
\begin{eqnarray*}
\langle [K^\bullet], [M^\bullet]\rangle&=&\langle [K^\bullet], [U_{\Ker(d_M)]}\rangle\langle [K^\bullet], [U_{\Im(d_M)}]\rangle\\
&=&\langle[K_{\Im(d_K)}], [U_{\Ker(d_M)}]\rangle\langle[K_{\Im(d_K)}], [U_{\Im(d_M)}]\rangle\\
&=&\langle [K_{\Im(d_K)}], [M^\bullet]\rangle.
\end{eqnarray*}
The proof of the other identity is similar.
\end{proof}

\section{Semi-derived Ringel-Hall algebras}\label{Section definition and structure}
Let $\varepsilon$ be an essentially small exact category, linear over a finite field.
Assume that $\varepsilon$ has finite morphism and extension spaces:
$$|\Hom_\varepsilon(A,B)|<\infty,\quad |\Ext^1_{\varepsilon}(A,B)|<\infty,\,\,\forall A,B\in\varepsilon.$$

Given objects $A,B,C\in\varepsilon$, define $\Ext_\varepsilon^1(A,C)_B\subseteq \Ext_\varepsilon^1(A,C)$ as the subset parameterizing extensions whose middle term is isomorphic to $B$. We define the Ringel-Hall algebra $\ch(\varepsilon)$ to be the $\Q$-vector space whose basis is parameterized by the isomorphism classes $[A]$ of objects $A$ of $\varepsilon$, with the multiplication
defined by
$$[A]\diamond [C]=\sum_{[B]\in \Iso(\varepsilon)}\frac{|\Ext^1_{\varepsilon}(A,C)_B|}{|\Hom_{\varepsilon}(A,C)|}[B].$$
It is well known that
the algebra $\ch(\varepsilon)$ is associative and unital. And the unit of the Ringel-Hall algebra is $[0]$, where $0$ is the zero object of $\varepsilon$, see \cite{R0} and also \cite{Sch2,P,Br} for details.
\subsection{Semi-derived Ringel-Hall algebras}
Let $\ch(\cc_{\Z/t}(\ca))$ be the Ringel-Hall algebra of $\cc_{\Z/t}(\ca)$, i.e., the multiplication is given by the following sum
$$[M_1^\bullet]\diamond[M_2^\bullet]=\sum_{[M^\bullet]\in\Iso(\cc_{\Z/t}(\ca))}\frac{|\Ext_{\cc_{\Z/t}(\ca)}^1(M_1^\bullet, M_2^\bullet)_{M^\bullet}|}{|\Hom_{\cc_{\Z/t}(\ca)}(M_1^\bullet, M_2^\bullet)|}[M^\bullet]$$
for any two objects $M_1^\bullet, M_2^\bullet\in\cc_{\Z/t}(\ca)$, one can see \cite{R0,R2,Sch1,Sch2} for more details. Let $J_{\Z/t}$ be the ideal of $\ch(\cc_{\Z/t}(\ca))$ generated by all the differences
$$[M_1^\bullet]-[M_2^{\bullet}],$$
if $M_1^\bullet$ and $M_2^{\bullet}$ satisfy $H^i(M_1^\bullet)\cong H^i(M_2^{\bullet})$ and $\widehat{\Im(d^i_{M_1^{\bullet}})}=\widehat{\Im(d^i_{M_2^{\bullet}})}$ for all $i\in\mathbb{Z}/t$. Due to a proof similar to that of Lemma 3.13 in \cite{LuP21}, it is easy to see that all the differences $[M_1^\bullet]-[M_2^{\bullet}]$ are $K_0(\cc_{\Z/t}(\ca))$-homogenous elements, so the quotient algebra $\ch(\cc_{\Z/t}(\ca))/J_{\Z/t}$ is also a $K_0(\cc_{\Z/t}(\ca))$-graded algebra and the induced multiplication is also denoted by $\diamond$.

Similar to \cite[Lemma 3.13]{LuP21}, we have the following lemma whose proof is omitted.
\begin{lemma}\label{lemma equivalences of some exact sequences}
Let $\cb$ be an abelian category. Given objects $$U^\bullet=(U^i,d_{U^\bullet}^i)_{i\in\Z/t}, V^\bullet=(V^i,d_{V^\bullet}^i)_{i\in\Z/t}, W^\bullet=(W^i,d_{W^\bullet}^i)_{i\in\Z/t}\in\cc_{\Z/t}(\cb),$$
if there is a short exact sequence $0\rightarrow U^\bullet\xrightarrow{h_1^\bullet} V^\bullet\xrightarrow{h_2^\bullet} W^\bullet\rightarrow 0$, then the following statements are equivalent.

(i) $0\rightarrow \Im(d_{U^\bullet}^i)\xrightarrow{t_1^i} \Im(d_{V^\bullet}^i)\xrightarrow{t_2^i} \Im(d_{W^\bullet}^i)\rightarrow 0$ is short exact for each $i\in\mathbb{Z}/t$.

(ii) $0\rightarrow \Ker(d_{U^\bullet}^i)\xrightarrow{s_1^i} \Ker(d_{V^\bullet}^i)\xrightarrow{s_2^i} \Ker(d_{W^\bullet}^i)\rightarrow 0$ is short exact for each $i\in\mathbb{Z}/t$.

(iii) $0\rightarrow H^i(U^\bullet)\xrightarrow{r_1^i} H^i(V^\bullet)\xrightarrow{r_2^i} H^i(W^\bullet)\rightarrow 0$ is short exact for each $i\in\mathbb{Z}/t$.

\noindent Here the morphisms are induced by $h_1^\bullet$ and $h_2^\bullet$, and $H^i(M^\bullet)(i\in\mathbb{Z}/t)$ denote the homologies of a complex $M^\bullet\in\cc_{\mathbb{Z}/t}(\cb)$.

In particular, if $U^\bullet$ or $W^\bullet$ is acyclic,
then for each $i\in\mathbb{Z}/t$
$$\widehat{\Im(d_{V^\bullet}^i)}=\widehat{\Im(d_{U^\bullet}^i)}+ \widehat{\Im(d_{W^\bullet}^i)}$$
in $K_0(\cb)$.
\end{lemma}

The following lemma is an easy consequence of Proposition \ref{proposition euler form of acyclic and any complex} and Lemma \ref{lemma equivalences of some exact sequences}.
\begin{lemma}\label{lemma multiplcation in quotient algebra}
For any $K^\bullet\in\cc_{\Z/t}^{ac}(\ca)$ and $M^\bullet\in\cc_{\Z/t}(\ca)$, then we have
\begin{eqnarray*}
[M^\bullet]\diamond [K^\bullet]=\frac{1}{\langle [M^\bullet],[K^\bullet]\rangle}[M^\bullet\oplus K^\bullet],&&[K^\bullet]\diamond[M^\bullet]=\frac{1}{\langle [K^\bullet],[M^\bullet]\rangle}[K^\bullet\oplus M^\bullet]
\end{eqnarray*}
in $\ch(\cc_{\Z/t}(\ca))/J_{\Z/t}$. In particular, for any $K_1^\bullet, K_2^\bullet\in\cc_{\Z/t}^{ac}(\ca)$, we have
\begin{eqnarray*}
&&[K_1^\bullet]\diamond [K_2^\bullet]=\frac{1}{\langle [K_1^\bullet], [K_2^\bullet]\rangle}[K_1^\bullet\oplus K_2^\bullet]
\end{eqnarray*}
in $\ch(\cc_{\Z/t}(\ca))/J_{\Z/t}$.
\end{lemma}

Similarly, we also set $S_{\Z/t}$ to be the subset of $\ch(\cc_{\Z/t}(\ca))/J_{\Z/t}$ formed by all $r[K^\bullet]$ with $r\in\mathbb{Q}^{\times}$ and $K^\bullet\in\cc_{\Z/t}^{ac}(\ca)$. It is routine to check that $S_{\Z/t}$ is a multiplicatively closed subset with the identity $[0]\in S_{\Z/t}$. And one can also check that the localization of $\ch(\cc_{\Z/t}(\ca))/J_{\Z/t}$ with respect to $S_{\Z/t}$ exists, which is denoted by $(\ch(\cc_{\Z/t}(\ca))/J_{\Z/t})[S_{\Z/t}^{-1}]$ and the multiplication is also denoted by $\diamond$.

\begin{definition}\label{def semi-derived hall algebra}
$(\ch(\cc_{\Z/t}(\ca))/J_{\Z/t})[S_{\Z/t}^{-1}]$ is called the ($\mathbb{Z}/t$-graded) semi-derived Ringel-Hall algebra of $\ca$, and denoted by $\cs\cd\ch_{\Z/t}(\ca)$.
\end{definition}

Let $I_{\Z/t}$ be the ideal of $\ch(\cc_{\Z/t}(\ca))$ generated by all the differences $[L^\bullet]-[K^\bullet\oplus M^\bullet]$, if there is a short exact sequence $K^\bullet\rightarrowtail L^\bullet\twoheadrightarrow M^\bullet$ in $\cc_{\Z/t}(\ca)$ with $K^\bullet$ acyclic. Obviously the ideal $I_{\Z/t}$ is contained in $J_{\Z/t}$. And it is the ideal $I_{\Z/t}$ that we used in our previous work \cite{LuP21} and \cite{LinP} when $t=2$ and $t=\infty$.
\begin{remark}
In case $t=2$ and $t=\infty$, the semi-derived Ringel-Hall algebra $\cs\cd\ch_{\Z/t}(\ca)$ is identified with $$(\ch(\cc_{\Z/t}(\ca))/I_{\Z/t})[S_{\Z/t}^{'-1}],$$ where $S_{\Z/t}'$ is the subset of $\ch(\cc_{\Z/t}(\ca))/I_{\Z/t}$ formed by all $r[K^\bullet]$ with $r\in\mathbb{Q}^{\times}$ and $K^\bullet\in\cc_{\Z/t}^{ac}(\ca)$. In fact we have a natural epimorphism from $(\ch(\cc_{\Z/t}(\ca))/I_{\Z/t})[S_{\Z/t}^{'-1}]$ to $\cs\cd\ch_{\Z/t}(\ca)$ since $I_{\Z/t}\subseteq J_{\Z/t}$. Note that a basis of $(\ch(\cc_{\Z/t}(\ca))/I_{\Z/t})[S_{\Z/t}^{'-1}]$ in \cite[Theorem 3.25]{LuP21} is the same as the one in Theorem \ref{theorem another basis in SDH(A)} of this paper (except for different symbols of stalk complexes) and the natural epimorphism sends the basis of $(\ch(\cc_{\Z/t}(\ca))/I_{\Z/t})[S_{\Z/t}^{'-1}]$ onto the basis of $\cs\cd\ch_{\Z/t}(\ca)$. So $(\ch(\cc_{\Z/t}(\ca))/I_{\Z/t})[S_{\Z/t}^{'-1}]$ is naturally isomorphic to $\cs\cd\ch_{\Z/t}(\ca)$. Similarly, for any $t\geq 1$ one has that $(\ch(\cc_{\Z/t}(\ca))/I_{\Z/t})[S_{\Z/t}^{'-1}]
\cong\cs\cd\ch_{\Z/t}(\ca)$. Our current definition of $\cs\cd\ch_{\Z/t}(\ca)$ is an improved version, which is more convenient for calculations.
\end{remark}

\begin{lemma}\label{lemma hall multiplicatoin of acyclic complexes}
For any $K^\bullet, K_1^\bullet, K_2^\bullet\in\cc_{\Z/t}^{ac}(\ca)$ and $M^\bullet\in\cc_{\Z/t}(\ca)$, we have
\begin{eqnarray*}
&&[K^\bullet]\diamond [M^\bullet]=\frac{1}{\langle [K^\bullet],[M^\bullet]\rangle}[K^\bullet\oplus M^\bullet]=\frac{\langle [M^\bullet],[K^\bullet]\rangle}{\langle [K^\bullet],[M^\bullet]\rangle}[M^\bullet]\diamond [K^\bullet]\\
&&[K_1^\bullet]^{-1}\diamond [K_2^\bullet]^{-1}=\langle [K_2^\bullet],[K_1^\bullet]\rangle [K_1^\bullet\oplus K_2^\bullet]^{-1},\\
&&[K_1^\bullet]^{-1}\diamond [K_2^\bullet]=\frac{\langle [K_1^\bullet], [K_2^\bullet]\rangle}{\langle [K_2^\bullet],[K_1^\bullet]\rangle}[K_2^\bullet]\diamond [K_1^\bullet]^{-1}
\end{eqnarray*}
in $\cs\cd\ch_{\Z/t}(\ca)$.
\end{lemma}

\subsection{A basis of $\cs\cd\ch_{\Z/t}(\ca)$}\label{The basis of the semi-derived Ringel-Hall algebras}
Consider the set $\Iso(\cc_{\Z/t}^{ac}(\ca))$ of isomorphism classes $[K^\bullet]$ of acyclic $\mathbb{Z}/t$-graded complexes and its quotient by the following set of relations:
$$\{[K_1^\bullet]=[K_2^\bullet]|\widehat{\Im(d^i_{K_1^\bullet})}=\widehat{\Im(d^i_{K_2^\bullet})}, \forall~i\in\Z/t\}.$$
If we endow $\Iso(\cc_{\Z/t}^{ac}(\ca))$ with the addition given by direct sums, this quotient gives the {\it Grothendieck monoid} denoted by $\cm_{\Z/t}^{ac}(\ca)$. And define the {\it quantum affine space} $\mathbb{A}_{\Z/t}^{ac}({\ca})$ of $\mathbb{Z}/t$-graded acyclic complexes to be the $\mathbb{Q}$-monoid algebra of $\cm_{\Z/t}^{ac}(\ca)$, with the multiplication twisted by the inverse of Euler form, i.e., the product of classes of acyclic complexes $K_1^\bullet,K_2^\bullet\in\cc^{ac}_{\Z/t}({\ca})$ is defined as follows:
$$[K_1^\bullet]\diamond[K_2^\bullet]:=\frac{1}{\langle K_1^\bullet,K_2^\bullet\rangle}[K_1^\bullet\oplus K_2^\bullet].$$
Clearly, $[K_1^\bullet\oplus K_2^\bullet]=[K_1^{'\bullet}\oplus K_2^{'\bullet}]$ if $[K_1^\bullet]=[K_1^{'\bullet}]$, $[K_2^\bullet]=[K_2^{'\bullet}]\in\cm_{\Z/t}^{ac}(\ca)$. As a result, $$[K_1^\bullet]\diamond[K_2^\bullet]=[K_1^{'\bullet}]\diamond[K_2^{'\bullet}]$$
follows from Proposition \ref{proposition euler form of acyclic and any complex}.

Define the {\it quantum torus} $\mathbb{T}_{\Z/t}^{ac}(\ca)$ of $\mathbb{Z}/t$-graded acyclic complexes to be the $\mathbb{Q}$-group algebra of $K_0(\mathcal{C}_{\Z/t}^{ac}(\ca))$ with the multiplication twisted by the Euler form as above. Note that $\mathbb{T}_{\Z/t}^{ac}(\ca)$ is the localization of $\mathbb{A}_{\Z/t}^{ac}({\ca})$ with respect to the set formed by all the classes of acyclic complexes and then it is generated by the classes $[K^\bullet]$ and their inverses $[K^\bullet]^{-1}$.

Note that $\cs\cd\ch_{\Z/t}(\ca)$ is a $\T_{\Z/t}^{ac}(\ca)$-bimodule with the bimodule structure induced by the Hall product.

We set $G_{\Z/t}$ to be the following set
$$\{(\alpha^\bullet,[A^\bullet])=(\alpha^0,\alpha^1,\cdots,\alpha^{t-1},[A^0],[A^1],\cdots,[A^{t-1}])\,|\, \alpha^i\in K_0(\ca), [A^i]\in \Iso(\ca), i\in\mathbb{Z}/t\}.$$
In the following, for any $(\alpha^\bullet,[A^\bullet])\in G_{\Z/t}$, we denote by $\cl_{(\alpha^\bullet,[A^\bullet])}$ the set of isomorphism class $[M^\bullet]$ of the $\mathbb{Z}/t$-graded complexes $$M^\bullet:= M^0\xrightarrow{d_{M^\bullet}^0}M^1\xrightarrow{d_{M^\bullet}^1}\cdots \xrightarrow{d_{M^\bullet}^{t-2}}M^{t-1}\xrightarrow{d_{M^\bullet}^{t-1}}M^0$$ satisfying $H^i(M^\bullet)\cong A^i$ and $\widehat{\Im(d_{M^\bullet}^i)}=\alpha^{i+1}$ for all $i\in\Z/t$.
So one can easily get that $\ch(\cc_{\Z/t}(\ca))$ is a $G_{\Z/t}$-graded linear space, i.e.,
$$\ch(\cc_{\Z/t}(\ca))=\bigoplus_{(\alpha^\bullet,[A^\bullet])\in G_{\Z/t}} (\bigoplus_{[M^\bullet]\in\cl_{(\alpha^\bullet,[A^\bullet])}}\Q[M^\bullet]).$$

We can similarly define $J'_{\Z/t}$ to be the linear subspace of $\ch(\cc_{\Z/t}(\ca))$ spanned by all the differences
$$[M_1^\bullet]-[M_2^{\bullet}]$$
satisfying $H^i(M_1^\bullet)\cong H^i(M_2^{\bullet})$ and $\widehat{\Im(d^i_{M_1^{\bullet}})}=\widehat{\Im(d^i_{M_2^{\bullet}})}$ for all $i\in\mathbb{Z}/t$. Clearly, $J_{\Z/t}'$ is a homogeneous subspace, so $\ch(\cc_{\Z/t}(\ca))/J_{\Z/t}'$ is also a $G_{\Z/t}$-graded linear space.

For any $A,B\in \ca$ with $\widehat{A}=\widehat{B}\in K_0(\ca)$ and $i\in\mathbb{Z}/t$, we have $[K_{A, i}]=[K_{B, i}]$ in both $\ch(\cc_{\Z/t}(\ca))/J_{\Z/t}$ and $\ch(\cc_{\Z/t}(\ca))/J'_{\Z/t}$ thus denoted by $K_{\widehat{A}, i}$. By definition the following lemma is clear.

\begin{lemma}\label{lemma acyclic in 1-periodic and iso of 1 acyclic of cone and K_M}
For any $M^\bullet=(M^i, d^i_{M^\bullet})_{i\in\mathbb{Z}/t}$, we have
$$[M^\bullet]=\left[\left(\bigoplus_{i\in\mathbb{Z}/t} K_{\widehat{\Im(d^i_{M^\bullet})},i+1}\right)\bigoplus\left(\bigoplus _{i\in\mathbb{Z}/t}U_{H^i(M^\bullet),i}\right)\right]$$
in $\ch(\cc_{\Z/t}(\ca))/J_{\Z/t}$ and $\ch(\cc_{\Z/t}(\ca))/J_{\Z/t}'$.
\end{lemma}

In addition, one can obtain a basis of $\ch(\cc_{\Z/t}(\ca))/J_{\Z/t}'$ in a way similar to \cite[Proposition 3.15]{LuP21}.
\begin{lemma}\label{lemma basis quotient space}
The linear space $\ch(\cc_{\Z/t}(\ca))/J_{\Z/t}'$ has a basis consisting of the following isomorphism classes
\[[K^\bullet\oplus U^\bullet], \]
for all $(\alpha^\bullet,[A^\bullet])\in G_{\Z/t}$, where $K^\bullet=\bigoplus_{i\in\Z/t}K_{\alpha^i,i}\in\cl_{(\alpha^\bullet,[0^\bullet])}$, $U^\bullet=\bigoplus_{i\in\Z/t}U_{A^i,i}\in\cl_{(0^\bullet,[A^\bullet])}$ and $\alpha^i\in K_0^+(\ca)(i\in\Z/t)$, here $K_{0}^+(\ca)$ denotes the positive cone in the Grothendieck group consisting of the classes of objects of $\ca$.
\end{lemma}

\begin{proposition}\label{proposition iso between quotient ideal and linear space}
$\ch(\cc_{\Z/t}(\ca))/J_{\Z/t}$ is isomorphic to $\ch(\cc_{\Z/t}(\ca))/J_{\Z/t}'$ as linear spaces.
\end{proposition}

\begin{proof}
Obviously $J_{\Z/t}'$ is contained in $J_{\Z/t}$.

On the other hand, for any $[M^\bullet]\in\cl_{(\alpha^\bullet,[A^\bullet])}$ it is easy to see that
$$q([M^\bullet])=[K^\bullet\oplus U^\bullet],$$
where $q: \ch(\cc_{\Z/t}(\ca))\rightarrow \ch(\cc_{\Z/t}(\ca))/J_{\Z/t}'$ is the natural projection, $K^\bullet=\bigoplus_{i\in\Z/t}K_{\alpha^i,i}\in\cl_{(\alpha^\bullet,[0^\bullet])}$ and $U^\bullet=\bigoplus_{i\in\Z/t}U_{A^i,i}\in\cl_{(0^\bullet,[A^\bullet])}$.

Let $[L^\bullet],[M_1^\bullet]$ and $[M_2^{\bullet}]$ be objects of $\Iso(\cc_{\Z/t}(\ca))$ satisfying $[M_1^{\bullet}], [M_2^{\bullet}]\in\cl_{(\alpha^\bullet,[A^\bullet])}$, then we claim that
$$q\left([L^\bullet]\diamond([M_1^\bullet]-[M_2^{\bullet}])\right)=0\ \text{and}\ q\left(([M_1^\bullet]-[M_2^{\bullet}])\diamond[L^\bullet]\right)=0.$$
We just show the first identity while the other one can be obtained dually. It is not hard to see that there is a complex $X^\bullet$ and acyclic complexes $K_1^\bullet,K_2^{\bullet}$ satisfying the following two short exact sequences
$$
0\rightarrow K_1^\bullet\rightarrow X^\bullet\rightarrow M_1^\bullet\rightarrow 0$$
and$$0\rightarrow K_2^{\bullet}\rightarrow X^\bullet\rightarrow M_2^{\bullet}\rightarrow 0.
$$
Without loss of generality, we can assume that $[X^\bullet]\in\cl_{(\beta^\bullet,[A^\bullet])}$, and then by Lemma \ref{lemma equivalences of some exact sequences} we know that
$$\widehat{\Im(d^i_{X^\bullet})}=\widehat{\Im(d^i_{K_1^\bullet})}+\widehat{\Im(d^i_{M_1^\bullet})}=\widehat{\Im(d^i_{K_2^\bullet})}+\widehat{\Im(d^i_{M_2^\bullet})}$$
for each $i\in\Z/t$. Thus, $[K_1^\bullet], [K_2^\bullet]\in\cl_{(\beta^\bullet-\alpha^\bullet,[0^\bullet])}$.

Owing to the following two identities
\begin{eqnarray*}
q([L^\bullet]\diamond[M_1^\bullet])&=&q\left(\sum_{[V_1^\bullet]\in\Iso(\cc_{\Z/t}(\ca))}\frac{|\Ext^1_{\cc_{\Z/t}(\ca)}(L^\bullet,M_1^\bullet)_{V_1^\bullet}|}{|\Hom_{\cc_{\Z/t}(\ca)}(L^\bullet,M_1^\bullet)|}[V_1^\bullet]\right)\\
&=&\sum_{(\gamma^\bullet,[C^\bullet])\in G_{\Z/t}}\left(\sum_{[V_1^\bullet]\in\cl_{(\gamma^\bullet,[C^\bullet])}}\frac{|\Ext^1_{\cc_{\Z/t}(\ca)}(L^\bullet,M_1^\bullet)_{V_1^\bullet}|}{|\Hom_{\cc_{\Z/t}(\ca)}(L^\bullet,M_1^\bullet)|}\right)q([V_1^\bullet])
\end{eqnarray*}
and
\begin{eqnarray*}
q([L^\bullet]\diamond[M_2^\bullet])&=&q\left(\sum_{[V_2^\bullet]\in\Iso(\cc_{\Z/t}(\ca))}\frac{|\Ext^1_{\cc_{\Z/t}(\ca)}(L^\bullet,M_2^\bullet)_{V_2^\bullet}|}{|\Hom_{\cc_{\Z/t}(\ca)}(L^\bullet,M_2^\bullet)|}[V_2^\bullet]\right)\\
&=&\sum_{(\gamma^\bullet,[C^\bullet])\in G_{\Z/t}}\left(\sum_{[V_2^\bullet]\in\cl_{(\gamma^\bullet,[C^\bullet])}}\frac{|\Ext^1_{\cc_{\Z/t}(\ca)}(L^\bullet,M_2^\bullet)_{V_2^\bullet}|}{|\Hom_{\cc_{\Z/t}(\ca)}(L^\bullet,M_2^\bullet)|}\right)q([V_2^\bullet]),
\end{eqnarray*}
one can easily see that
$$q([V_1^\bullet])=q( [V_2^\bullet])$$
for any $[V_1^\bullet], [V_2^\bullet]\in\cl_{(\gamma^\bullet,[C^\bullet])}$.
We just need to prove that
$$\sum_{[V_1^\bullet]\in\cl_{(\gamma^\bullet,[C^\bullet])}}\frac{|\Ext^1_{\cc_{\Z/t}(\ca)}(L^\bullet,M_1^\bullet)_{V_1^\bullet}|}{|\Hom_{\cc_{\Z/t}(\ca)}(L^\bullet,M_1^\bullet)|}
=\sum_{[V_2^\bullet]\in\cl_{(\gamma^\bullet,[C^\bullet])}}\frac{|\Ext^1_{\cc_{\Z/t}(\ca)}(L^\bullet,M_2^\bullet)_{V_2^\bullet}|}{|\Hom_{\cc_{\Z/t}(\ca)}(L^\bullet,M_2^\bullet)|}$$ for any $(\gamma^\bullet,[C^\bullet])\in G_{\Z/t}$.
Obviously we have the following two exact sequences

$0\rightarrow \Hom_{\cc_{\Z/t}(\ca)}(L^\bullet, K_1^\bullet)\rightarrow \Hom_{\cc_{\Z/t}(\ca)}(L^\bullet, X^\bullet)\rightarrow\Hom_{\cc_{\Z/t}(\ca)}(L^\bullet, M_1^\bullet)\rightarrow\Ext^1_{\cc_{\Z/t}(\ca)}(L^\bullet, K_1^\bullet)\rightarrow\Ext^1_{\cc_{\Z/t}(\ca)}(L^\bullet, X^\bullet)\xrightarrow{\phi}\Ext^1_{\cc_{\Z/t}(\ca)}(L^\bullet, M_1^\bullet)\rightarrow 0$ and

$0\rightarrow \Hom_{\cc_{\Z/t}(\ca)}(L^\bullet, K_2^\bullet)\rightarrow \Hom_{\cc_{\Z/t}(\ca)}(L^\bullet, X^\bullet)\rightarrow\Hom_{\cc_{\Z/t}(\ca)}(L^\bullet, M_2^\bullet)\rightarrow\Ext^1_{\cc_{\Z/t}(\ca)}(L^\bullet, K_2^\bullet)\rightarrow\Ext^1_{\cc_{\Z/t}(\ca)}(L^\bullet, X^\bullet)\xrightarrow{\omega}\Ext^1_{\cc_{\Z/t}(\ca)}(L^\bullet, M_2^\bullet)\rightarrow 0$.\\
Thus by the similar proof of Lemma 3.16 in \cite{LuP21}, we get the following two epimorphisms
$$\Phi: \coprod_{[W^\bullet]\in\cl_{(\beta^\bullet-\alpha^\bullet+\gamma^\bullet,[C^\bullet])}}\Ext^1_{\cc_{\Z/t}(\ca)}(L^\bullet, X^\bullet)_{W^\bullet}\rightarrow\coprod_{[V_1^\bullet]\in\cl_{(\gamma^\bullet,[C^\bullet])}}\Ext^1_{\cc_{\Z/t}(\ca)}(L^\bullet, M_1^\bullet)_{V_1^\bullet}$$
and
$$\Omega: \coprod_{[W^\bullet]\in\cl_{(\beta^\bullet-\alpha^\bullet+\gamma^\bullet,[C^\bullet])}}\Ext^1_{\cc_{\Z/t}(\ca)}(L^\bullet, X^\bullet)_{W^\bullet}\rightarrow\coprod_{[V_2^\bullet]\in\cl_{(\gamma^\bullet,[C^\bullet])}}\Ext^1_{\cc_{\Z/t}(\ca)}(L^\bullet, M_2^\bullet)_{V_2^\bullet},$$
with $$|\Phi^{-1}(0)|=|\Ker(\phi)|=\frac{|\Hom_{\cc_{\Z/t}(\ca)}(L^\bullet, X^\bullet)||\Ext^1_{\cc_{\Z/t}(\ca)}(L^\bullet, K_1^\bullet)|}{|\Hom_{\cc_{\Z/t}(\ca)}(L^\bullet, K_1^\bullet)||\Hom_{\cc_{\Z/t}(\ca)}(L^\bullet, M_1^\bullet)|}$$ and $$|\Omega^{-1}(0)|=|\Ker(\omega)|=\frac{|\Hom_{\cc_{\Z/t}(\ca)}(L^\bullet, X^\bullet)||\Ext^1_{\cc_{\Z/t}(\ca)}(L^\bullet, K_2^\bullet)|}{|\Hom_{\cc_{\Z/t}(\ca)}(L^\bullet, K_2^\bullet)||\Hom_{\cc_{\Z/t}(\ca)}(L^\bullet, M_2^\bullet)|}.$$
So by Proposition \ref{proposition euler form of acyclic and any complex} we have that
\begin{eqnarray*}
&&\sum_{[V_1^\bullet]\in\cl_{(\gamma^\bullet,[C^\bullet])}}\frac{|\Ext^1_{\cc_{\Z/t}(\ca)}(L^\bullet,M_1^\bullet)_{V_1^\bullet}|}{|\Hom_{\cc_{\Z/t}(\ca)}(L^\bullet,M_1^\bullet)|}\\
&=&\sum_{[W^\bullet]\in\cl_{(\beta^\bullet-\alpha^\bullet+\gamma^\bullet,[C^\bullet])}}\frac{|\Ext^1_{\cc_{\Z/t}(\ca)}(L^\bullet,X^\bullet)_{W^\bullet}|}{|\Phi^{-1}(0)||\Hom_{\cc_{\Z/t}(\ca)}(L^\bullet,M_1^\bullet)|}\\
&=&\\
&=&\sum_{[W^\bullet]\in\cl_{(\beta^\bullet-\alpha^\bullet+\gamma^\bullet,[C^\bullet])}}\frac{|\Ext^1_{\cc_{\Z/t}(\ca)}(L^\bullet,X^\bullet)_{W^\bullet}|}{|\Hom_{\cc_{\Z/t}(\ca)}(L^\bullet,X^\bullet)|}\langle L^\bullet, K_1^\bullet\rangle\\
&=&\sum_{[W^\bullet]\in\cl_{(\beta^\bullet-\alpha^\bullet+\gamma^\bullet,[C^\bullet])}}\frac{|\Ext^1_{\cc_{\Z/t}(\ca)}(L^\bullet,X^\bullet)_{W^\bullet}|}{|\Hom_{\cc_{\Z/t}(\ca)}(L^\bullet,X^\bullet)|}\langle L^\bullet, K_2^\bullet\rangle\\
&=&\sum_{[V_2^\bullet]\in\cl_{(\gamma^\bullet,[C^\bullet])}}\frac{|\Ext^1_{\cc_{\Z/t}(\ca)}(L^\bullet,M_2^\bullet)_{V_2^\bullet}|}{|\Hom_{\cc_{\Z/t}(\ca)}(L^\bullet,M_2^\bullet)|}.
\end{eqnarray*}

Therefore one can obtain the induced linear map $\bar{q}: \ch(\cc_{\Z/t}(\ca))/J_{\Z/t}\rightarrow \ch(\cc_{\Z/t}(\ca))/J_{\Z/t}'$ which implies that $J_{\Z/t}\subseteq J_{\Z/t}'$. This completes the proof.
\end{proof}

And we get the following results immediately from Lemma \ref{lemma acyclic in 1-periodic and iso of 1 acyclic of cone and K_M}, Proposition \ref{proposition iso between quotient ideal and linear space} and Proposition \ref{prop euler form of acycli complexes}.

\begin{proposition}\label{proposition decomposition in semi-derived}
Let $M^\bullet= M^0\xrightarrow{d^0_{M^\bullet}}M^1\xrightarrow{d^1_{M^\bullet}}\cdots \xrightarrow{d^{t-2}_{M^\bullet}}M^{t-1}\xrightarrow{d^{t-1}_{M^\bullet}}M^0$ be a $\mathbb{Z}/t$-graded complex. Then
we have
\begin{eqnarray*}
[M]&=&\langle\widehat{\Im(d^{t-1}_{M^\bullet})}, \widehat{\Im(d^{t-2}_{M^\bullet})}\rangle\prod_{i=0}^{t-1}\langle\widehat{\Im(d^i_{M^\bullet})}, \widehat{H^i(M^\bullet)}\rangle\\
&&[K_{\widehat{\Im(d^{t-1}_{M^\bullet})}, 0}]\diamond[K_{\widehat{\Im(d^{0}_{M^\bullet})}, 1}]\diamond\cdots\diamond[K_{\widehat{\Im(d^{t-2}_{M^\bullet})}, t-1}]\diamond[\bigoplus\limits_{i=0}^{t-1}U_{H^i(M^\bullet),i}]
\end{eqnarray*}
in $(\ch(\cc_{\Z/t}(\ca))/J_{\Z/t}$.
\end{proposition}

Induced by Lemma \ref{lemma basis quotient space}, Proposition \ref{proposition iso between quotient ideal and linear space} and Proposition \ref{proposition decomposition in semi-derived}, one can obtain a basis of $\ch(\cc_{\Z/t}(\ca))/J_{\Z/t}$ as follows .

\begin{proposition}\label{proposition basis in quotient algebra}
The algebra $\ch(\cc_{\Z/t}(\ca))/J_{\Z/t}$ has a basis consisting of the following elements
\[[K_{\alpha^0, 0}]\diamond[K_{\alpha^1, 1}]\diamond\cdots\diamond[K_{\alpha^{t-1}, t-1}]\diamond[\bigoplus\limits_{i=0}^{t-1}U_{A^i,i}],\]
where $\alpha^i\in K_{0}^+(\ca)$, $[A^i]\in\Iso(\ca)$ for all $i\in\Z/t$.
\end{proposition}

\begin{corollary}\label{cor free module over quantum affine}
The algebra $\ch(\cc_{\Z/t}(\ca))/J_{\Z/t}$ is a free module over the quantum affine space $\mathbb{A}_{\Z/t}^{ac}(\ca)$, with a basis given by $$\{[\bigoplus\limits_{i=0}^{t-1}U_{A^i,i}]|[A^i]\in\Iso(\ca), i\in\mathbb{Z}/t\}.$$
\end{corollary}

If $\alpha=\widehat{A}-\widehat{B}$, set $$K_{\alpha, i}:=\frac{1}{\langle\alpha, \widehat{B}\rangle}[K_{\widehat{A}, i}]\diamond[K_{\widehat{B}, i}]^{-1}$$ in case $t>1$ and $$K_\alpha:=\frac{1}{\langle\alpha, \widehat{B}\rangle^2}[K_{\widehat{A}}]\diamond[K_{\widehat{B}}]^{-1}$$
in case $t=1$. Following the similar proof in \cite{LuP21}, we know that both $K_{\alpha, i}$ and $K_\alpha$ are well-defined, i.e., they are all independent of the expression of $\alpha$ in $K_0(\ca)$, and the identities in Proposition \ref{prop euler form of acycli complexes} are all valid for all $K_{\alpha, i}$ and $K_{\alpha}$.

\begin{theorem}\label{theorem basis in VS}
$\cs\cd\ch_{\Z/t}(\ca)$ has a basis consisting of the following elements
$$[K_{\alpha^0, 0}]\diamond[K_{\alpha^1, 1}]\diamond\cdots\diamond[K_{\alpha^{t-1}, t-1}]\diamond[\bigoplus\limits_{i=0}^{t-1}U_{A^i,i}],$$
where $\alpha^i\in K_0(\ca), [A^i]\in\Iso(\ca)$ for any $i\in\Z/t$.
\end{theorem}

\begin{corollary}\label{cor free module over torus}
$\cs\cd\ch_{\Z/t}(\ca)$ is a free module over the quantum torus $\mathbb{T}_{\Z/t}^{ac}(\ca)$, with a basis given by $$\{[\bigoplus\limits_{i=0}^{t-1}U_{A^i,i}]|[A^i]\in\Iso(\ca), i\in\mathbb{Z}/t\}.$$
\end{corollary}

The following is analogous to Theorem 3.25 in \cite{LuP21}, and we omit the proof.
\begin{theorem}\label{theorem another basis in SDH(A)}
$\cs\cd\ch_{\Z/t}(\ca)$ has a basis given by the following elements
$$[K_{\alpha^0, 0}]\diamond[K_{\alpha^1, 1}]\diamond\cdots\diamond[K_{\alpha^{t-1}, t-1}]\diamond[U_{A^{t-1},t-1}]\diamond\cdots\diamond[U_{A^{0},0}],$$
where $\alpha^i\in K_0(\ca), [A^i]\in\Iso(\ca)$ for any $i\in\Z/t$.
\end{theorem}
\subsection{Generators and defining relations of the semi-derived Ringel-Hall algebras}
For any $i\in\mathbb{Z}/t$, we also use $K_{\alpha,i}$ and $U_{A,i}$ to denote their isomorphism classes $[K_{\alpha,i}]$ and $[U_{A,i}]$ respectively if there is no risk of confusion.
Given four objects $A, B, M$ and $N$ of $\ca$, let $V(M, B, A, N)$ be the subset of $\Hom_\ca(M,B)\times\Hom_\ca(B,A)\times\Hom_\ca(A,N)$ consisting of exact sequences $0\rightarrow M \rightarrow B\rightarrow A\rightarrow N\rightarrow 0$. The set $V(M, B, A, N)$ is finite and we define a rational number $\gamma_{AB}^{MN}:=\frac{|V(M, B, A, N)|}{a_Aa_B}$, where $a_A=|\aut(A)|$ and $a_B=|\aut(B)|$.
\begin{proposition}\label{proposition generators 1-periodic}
The semi-derived Ringel-Hall algebra $\cs\cd\ch_{\Z/1}(\ca)$ is generated by the set $$\{U_A, K_{\alpha}|[A]\in\Iso(\ca), \alpha\in K_0(\ca)\},$$ with the following defining relations (\ref{relation 1 in MH1})-(\ref{relation 3 in MH1}).
\begin{eqnarray}
U_A\diamond U_B&=&\sum\limits_{C\in\Iso(\ca)}
\left(\sum\limits_{\tiny \begin{array}{c}[L^\bullet]\in\Iso(\cc_{\Z/1}(\ca))\\H(L^\bullet)\cong C\end{array}}\frac{|\Ext^1_{\cc_{\Z/1}(\ca)}(U_A, U_B)_{L^\bullet}|}{|\Hom_{\ca}(A, B)|}\right)\label{relation 1 in MH1}\\
&&\ \ \ \sqrt{\langle\widehat{A}+\widehat{B}-\widehat{C}, \widehat{C}\rangle} K_{\frac{1}{2}(\widehat{A}+\widehat{B}-\widehat{C})}\diamond U_C,\nonumber\\
U_A\diamond K_\alpha&=&\frac{\langle\alpha, \widehat{A}\rangle}{\langle\widehat{A},\alpha\rangle}K_\alpha\diamond U_A,\label{relation 2 in MH1}\\
K_\alpha\diamond K_\beta&=&\frac{1}{\langle\alpha,\beta\rangle^2}K_{\alpha+\beta}, \label{relation 3 in MH1}
\end{eqnarray}
where $[A], [B]\in\Iso(\ca)$ and $\alpha,\beta\in K_0(\ca)$.
\end{proposition}
\begin{proof}
For any $A, B\in\ca$ and any short exact sequence
$$0\rightarrow U_B \rightarrow L^\bullet\rightarrow U_A\rightarrow 0,$$
where $L^\bullet=(L,d_{L})$,
one can easily see that $\widehat{\Im(d_{L})}=\frac{1}{2}(\widehat{A}+\widehat{B}-\widehat{H(L^\bullet))}$. So according to Proposition \ref{proposition decomposition in semi-derived} we have
$$[L^\bullet]=\sqrt{\langle\widehat{A}+\widehat{B}-\widehat{H(L^\bullet)}, \widehat{H(L^\bullet)}\rangle}K_{\frac{1}{2}(\widehat{A}+\widehat{B}-\widehat{H(L^\bullet)})}\diamond U_{H(L^\bullet)}$$
in $\cs\cd\ch_{\Z/1}(\ca)$, thus
\begin{eqnarray*}
U_A\diamond U_B&=&\sum\limits_{C\in\Iso(\ca)}
\left(\sum\limits_{\tiny \begin{array}{c}[L^\bullet]\in\Iso(\cc_{\Z/1}(\ca))\\H(L^\bullet)\cong C\end{array}}\frac{|\Ext^1_{\cc_{\Z/1}(\ca)}(U_A, U_B)_{L^\bullet}|}{|\Hom_{\ca}(A, B)|}\right) \sqrt{\langle\widehat{A}+\widehat{B}-\widehat{C}, \widehat{C}\rangle}K_{\frac{1}{2}(\widehat{A}+\widehat{B}-\widehat{C})}\diamond U_C.
\end{eqnarray*}
And it is easy to obtain the remaining relations by the definition of the semi-derived Ringel-Hall algebras and Proposition \ref{prop euler form of acycli complexes}.
By the basis given in Theorem \ref{theorem basis in VS}, the relations (\ref{relation 1 in MH1})-(\ref{relation 3 in MH1}) are actually the defining relations of $\cs\cd\ch_{\Z/1}(\ca)$.
\end{proof}
And the next two propositions are both obvious by the similar proof in \cite{LinP}.
\begin{proposition}\label{proposition generators 3-periodic}
The semi-derived Ringel-Hall algebra $\cs\cd\ch_{\Z/3}(\ca)$ is generated by the set $$\{U_{A,i}, K_{\alpha,i}|[A]\in\Iso(\ca), \alpha\in K_0(\ca), i\in\mathbb{Z}/3\},$$ with the following defining relations (\ref{relation 1 in MH3})-(\ref{relation 6 in MH3}).
\begin{eqnarray}
U_{A,i}\diamond U_{B,i}&=&\sum\limits_{C\in\Iso(\ca)}\frac{|\Ext_{\ca}^{1}(A,B)_{C}|}{|\Hom_{\ca}(A,B)|}U_{C,i},\label{relation 1 in MH3}\\
K_{\alpha,i}\diamond U_{B,i}&=&\langle \widehat{B},\alpha\rangle U_{B,i}\diamond K_{\alpha,i},\\
K_{\alpha,i}\diamond K_{\beta,i}&=&\frac{1}{\langle\alpha,\beta\rangle}K_{\alpha+\beta,i},
\end{eqnarray}

\begin{eqnarray}
U_{B,i}\diamond U_{A,i+1}&=&\sum\limits_{M,N\in\Iso(\ca)}\gamma_{AB}^{MN}\frac{a_Aa_B}{a_Ma_N}\langle \widehat{B}-\widehat{M}, \widehat{M}\rangle K_{\widehat{B}-\widehat{M}, i+1}\diamond U_{N,i+1}\diamond U_{M, i},
\end{eqnarray}

\begin{eqnarray}
U_{B,i}\diamond K_{\alpha,i+1}&=&\langle\alpha, \widehat{B}\rangle K_{\alpha,i+1}\diamond U_{B,i},\\
K_{\alpha,i}\diamond U_{B,i+1}&=&U_{B,i+1}\diamond K_{\alpha,i},\\
K_{\beta,i}\diamond K_{\alpha,i+1}&=&\langle\alpha,\beta\rangle K_{\alpha,i+1}\diamond K_{\beta,i},\label{relation 6 in MH3}
\end{eqnarray}
where $[A], [B]\in\Iso(\ca)$, $\alpha,\beta\in K_0(\ca)$ and $i\in\mathbb{Z}/3$.
\end{proposition}
\begin{proposition}\label{proposition generators t-periodic}The semi-derived Ringel-Hall algebra $\cs\cd\ch_{\Z/t}(\ca)(t>3)$ is generated by the set $$\{U_{A,i}, K_{\alpha,i}|A\in\Iso(\ca), \alpha\in K_0(\ca), i\in\mathbb{Z}/t\},$$ with the defining relations (\ref{relation 1 in MH3})-(\ref{relation 10 in MHt}).
\begin{eqnarray}
U_{A,i}\diamond U_{B,j}&=&U_{B,j}\diamond U_{A,i},\label{relation 8 in MHt}\\
K_{\alpha,i}\diamond U_{B,j}&=&U_{B,j}\diamond K_{\alpha,i},\\
K_{\alpha,i}\diamond K_{\beta,j}&=&K_{\beta,j}\diamond K_{\alpha,i},\label{relation 10 in MHt}
\end{eqnarray}
where $[A], [B]\in\Iso(\ca)$, $\alpha,\beta\in K_0(\ca)$ and $i, j\in\mathbb{Z}/t$ with $j\neq i, i\pm1$.
\end{proposition}

\begin{remark}
Zhang described the Bridgeland's Hall algebra by the generators and generating relations in the case of $t>2$ in \cite{Zhang1}, however he did not show the correspondent results in the case of $t=1$ or $t=2$. To realize the universal $q$-Onsager algebra via $\imath$Hall algebra of the projective line Lu, Ruan and Wang defined the $\imath$Hall algebra in the case of $t=1$ in \cite{LRW20} which is by definition the twisted semi-derived Ringel-Hall algebra of the category of $1$-periodic complexes of coherent sheaves on the projective line, but they do not show its generators and generating relations. In the case of $t=2$, since a twisted version of the semi-derived Ringel-Hall algebra is isomorphic to the Drinfeld double of the twisted extended Ringel-Hall algebra \cite{LuP16,LuP21}, one can also describe the semi-derived Ringel-Hall algebra via its generators and generating relations.
\end{remark}

\section{Hall algebras of odd-periodic relative derived categories}\label{Hall algebras of odd-periodic relative derived categories}
Let $\ct$ be a $k$-additive triangulated category with the translation $T=[1]$ satisfying the following conditions:
\begin{itemize}
\item[(i)] $\dim_k\Hom_{\ct}(X, Y)<\infty$ for any two objects $X$ and $Y$;
\item[(ii)] $\End_{\ct}(X)$ is local for any indecomposable object $X$;
\item[(iii)] $T^t=[t]\cong 1_{\ct}$ for some positive odd integer $t$.
\end{itemize}
The category $\ct$ is called a $t$-periodic triangulated category. And the derived Hall algebra $\cd\ch(\ct)$ of the $t$-periodic triangulated category $\ct$ introduced by Xu and Chen in \cite{XC} is the $\mathbb{Q}(v)$-space with the basis $\{[X]|X\in\ct\}$ and the multiplication is defined by
$$[X][Y]=\sum_{[L]\in\Iso(\ct)}\frac{|\Ext^1_\ct(X, Y)_L|}{\sqrt{\prod\limits_{i=0}^{t-1}|\Hom_{\ct}(X[i], Y)|^{(-1)^i}}}[L],$$
where $\Ext^1_\ct(X, Y)_L$ is defined to be $\Hom_{\ct}(X, Y[1])_{L[1]}$ which denotes the subset of $\Hom_{\ct}(X, Y[1])$ consisting of morphisms $l: X\rightarrow Y[1]$ whose cone $\Cone(l)$ is isomorphic to $L[1]$. Here we use the structure constant slightly different from the original one, but it can be obtained from $\cd\ch(\ct)$ via substituting the basis $[X]$ by $\frac{1}{a_X\sqrt{\prod\limits_{i=0}^{t-1}|\Hom_{\ct}(X[i], X)|^{(-1)^i}}}[X]$ for all $[X]\in\Iso(\ct)$, and one can see \cite{SCX,XX2} for more details. Obviously we can define the derived Hall algebra of $\cd_{\Z/t}(\ca)$ in case $t$ is an odd integer, and we denote it by $\cd\ch_{\Z/t}(\ca)$.

For two tuples of integers $(a_i)$ and $(b_i)\in\mathbb{Z}^t$, we say $(a_i)<(b_i)$ if $a_i\leq b_i$ for any $i\in \mathbb{Z}/t$ and there exists at least one $j\in\mathbb{Z}/t$ satisfying $a_j<b_j$. Similar to the work of Guo and Peng in \cite{GP}, we can also define a partial order over $\Iso(\cd_{t}(\ca))$. By Proposition \ref{proposition t-periodic iso to the homology}, for any isomorphism classes $[M^\bullet]=[\bigoplus\limits_{i=0}^{t-1}U_{M^i,i}]$ and $[N^\bullet]=[\bigoplus\limits_{i=0}^{t-1}U_{N^i,i}]$ in $\Iso(\cd_{\Z/t}(\ca))$, we say
$$[M^\bullet]<[N^\bullet],$$
if $\sum_{i=0}^{t-1}\widehat{M^i}<\sum_{i=0}^{t-1}\widehat{N^i}$.

Following Theorem 3.4 in \cite{XC}, we can similarly describe the derived Hall algebra $\cd\ch_{\Z/t}(\ca)$ in terms of its generators and relations as follows by induction on the above partial order.

\begin{proposition}\label{proposition generators and relations 1-derived hall algebra}
$\cd\ch_{\Z/1}(\ca)$ is isomorphic to an associative and unital $\mathbb{Q}(v)$-algebra generated by the set $$\{Z_A|A\in\Iso(\ca)\},$$  with the defining relations (\ref{re1 in DH1})
\begin{eqnarray}
Z_AZ_B&=&\sum\limits_{C\in\Iso(\ca)}\frac{|\Ext^1_{\cd_{\Z/1}(\ca)}(Z_A,Z_B)_{Z_C}|}{\sqrt{|\Hom_{\cd_{\Z/1}(\ca)}(Z_A,Z_B)|}}Z_C,\label{re1 in DH1}
\end{eqnarray}
where $[A], [B]\in\Iso(\ca)$.
\end{proposition}

\begin{proposition}\label{proposition generators and relations 3-derived hall algebra}
$\cd\ch_{\Z/3}(\ca)$ is isomorphic to an associative and unital $\mathbb{Q}(v)$-algebra generated by the set $$\{Z_A^{[i]}|A\in\Iso(\ca),i\in\mathbb{Z}/3\},$$  with the defining relations as follows.
\begin{eqnarray}
Z_A^{[i]}Z_B^{[i]}&=&\sum\limits_{C\in \Iso(\ca)}\frac{|\Ext^1_{\ca}(A,B)_C|}{\sqrt{|\Hom_\ca(A,B)||\Ext_{\ca}^1(A,B)|}}Z_C^{[i]},\label{prop relation 1 in DH3}\\
Z_B^{[i]}Z_A^{[i+1]}&=&\sum\limits_{M, N\in\Iso(\ca)}\gamma_{AB}^{MN}\frac{a_Aa_B}{a_Ma_N}\frac{1}{\sqrt{\langle\widehat{B}, \widehat{A}\rangle\langle \widehat{N}, \widehat{M}\rangle}}Z_N^{[i+1]}Z_M^{[i]},\label{prop relation 2 in DH3}
\end{eqnarray}
where $[A], [B]\in\Iso(\ca)$ and $i\in\mathbb{Z}/3$.
\end{proposition}

\begin{proposition}\label{proposition generators and relations t-derived hall algebra}
$\cd\ch_{\Z/t}(\ca)(t>3)$ is isomorphic to an associative and unital $\mathbb{Q}(v)$-algebra generated by the set $$\{Z_A^{[i]}|A\in\Iso(\ca),i\in\mathbb{Z}/t\},$$  with the defining relations (\ref{prop relation 1 in DH3})-(\ref{relation 3 in DHt}).
\begin{eqnarray}
Z_A^{[i]}Z_B^{[j]}&=&\sqrt{(\widehat{A}, \widehat{B})^{(-1)^{j-i}}}Z_B^{[j]}Z_A^{[i]}\ \text{for any~} j=i+2,\cdots,t-1,\label{relation 3 in DHt}
\end{eqnarray}
where $[A], [B]\in\Iso(\ca)$ and $i,j\in\mathbb{Z}/t$.
\end{proposition}
\section{Main result}\label{main results}
This section is devoted to investigating a relationship between the semi-derived Ringel-Hall algebra $\cs\cd\ch_{\Z/t}(\ca)$ and the derived-Hall algebra $\cd\ch_{\Z/t}(\ca)$ in case $t$ is an odd integer.
\subsection{Extended semi-derived Ringel-Hall algebras}We first define the extended quantum torus $\mathbb{T}_{\Z/t,\frac{1}{2}}^{ac}(\ca)$ of acyclic $\mathbb{Z}/t$-graded complexes as in Section \ref{The basis of the semi-derived Ringel-Hall algebras}. Consider the set
$$\sqrt{\Iso(\cc_{\Z/t}^{ac}(\ca))}:=\{\sqrt{[K^\bullet]}|[K^\bullet]\in\Iso(\cc_{\Z/t}^{ac}(\ca))\}$$
and its quotient by the following set of relations
$$\{\sqrt{[K_1^\bullet]}=\sqrt{[K_2^\bullet]}|\widehat{\Im(d^i_{K_1^\bullet})}=\widehat{\Im(d^i_{K_2^\bullet})}, \forall~i\in\mathbb{Z}/t\}.$$
If we endow $\sqrt{\Iso(\cc_{\Z/t}^{ac}(\ca))}$ with the addition given by direct sums, i.e.,
$$\sqrt{[K_1^\bullet]}+\sqrt{[K_2^\bullet]}:=\sqrt{[K_1^\bullet\oplus K_2^\bullet]}$$
for any $\sqrt{[K_1^\bullet]},\sqrt{[K_2^\bullet]}\in\sqrt{\Iso(\cc_{\Z/t}^{ac}(\ca))}$, this quotient gives the extended monoid denoted by $\cm_{\Z/t,\frac{1}{2}}^{ac}(\ca)$. And define the extended quantum affine space $\mathbb{A}_{\Z/t,\frac{1}{2}}^{ac}({\ca})$ of acyclic $\mathbb{Z}/t$-graded complexes to be the $\mathbb{R}$-monoid algebra of $\cm_{\Z/t,\frac{1}{2}}^{ac}(\ca)$, with the multiplication twisted by the inverse of Euler form, i.e.,
$$\sqrt{[K_1^\bullet]}\diamond\sqrt{[K_2^\bullet]}:=\frac{1}{\sqrt[4]{\langle [K_1^\bullet],[K_2^\bullet]\rangle}}\sqrt{[K_1^\bullet\oplus K_2^\bullet]}$$
for any $K_1^\bullet,K_2^\bullet\in\cc^{ac}_{\Z/t}({\ca})$. It is clear that $\sqrt{[K_1^\bullet\oplus K_2^\bullet]}=\sqrt{[N_1^\bullet\oplus N_2^\bullet]}$, if $\sqrt{[K_1^\bullet]}=\sqrt{[N_1^\bullet]}$, $\sqrt{[K_2^\bullet]}=\sqrt{[N_2^\bullet]}\in\cm_{\Z/t,\frac{1}{2}}^{ac}(\ca)$. As a result, $$\sqrt{[K_1^\bullet]}\diamond\sqrt{[K_2^\bullet]}=\sqrt{[N_1^\bullet]}\diamond\sqrt{[N_2^\bullet]}$$
follows from Proposition \ref{proposition euler form of acyclic and any complex}.

Define the extended quantum torus $\mathbb{T}_{\Z/t,\frac{1}{2}}^{ac}(\ca)$ of acyclic $\mathbb{Z}/t$-graded complexes as the localization of $\mathbb{A}_{\Z/t,\frac{1}{2}}^{ac}({\ca})$ with respect to the set formed by all the classes of $\sqrt{[K^\bullet]}$ with $[K^\bullet]\in\Iso(\cc_{\Z/t}^{ac}(\ca))$.
If we set further $[K^\bullet]=\sqrt{[K^\bullet\oplus K^\bullet]}$, then we get the quantum torus $\mathbb{T}_{\Z/t}^{ac}(\ca)$.

\begin{lemma}\label{lemma hall multiplicatoin of acyclic complexes}
For any $K_1^\bullet, K_2^\bullet\in\cc_{\Z/t}^{ac}(\ca)$, we have
\begin{eqnarray*}
&&\sqrt{[K_1^\bullet]}^{-1}\diamond \sqrt{[K_2^\bullet]}^{-1}=\sqrt[4]{\langle [K_2^\bullet],[K_1^\bullet]\rangle} \sqrt{[K_1^\bullet\oplus K_2^\bullet]}^{-1},\\
&&\sqrt{[K_1^\bullet]}^{-1}\diamond \sqrt{[K_2^\bullet]}=\frac{\sqrt[4]{\langle [K_1^\bullet], [K_2^\bullet]\rangle}}{\sqrt[4]{\langle [K_2^\bullet],[K_1^\bullet]\rangle}}\sqrt{[K_2^\bullet]}\diamond \sqrt{[K_1^\bullet]}^{-1}
\end{eqnarray*}
in $\mathbb{T}_{\Z/t,\frac{1}{2}}^{ac}(\ca)$.
\end{lemma}

Similarly, set $$\sqrt{K_{\alpha, i}}:=\frac{1}{\sqrt[4]{\langle\alpha, \widehat{B}\rangle}}\sqrt{[K_{\widehat{A}, i}]}\diamond\sqrt{[K_{\widehat{B}, i}]}^{-1}$$ in case $t>1$ and $$\sqrt{K_\alpha}:=\frac{1}{\sqrt{\langle\alpha, \widehat{B}\rangle}}\sqrt{[K_{\widehat{A}}]}\diamond\sqrt{[K_{\widehat{B}}]}^{-1}$$
in case $t=1$ for any $\alpha=\widehat{A}-\widehat{B}\in K_0(\ca)$. For any $\alpha,\beta\in K_0(\ca)$, it is routine to check that
$$\sqrt{K_{\alpha, i}}\diamond\sqrt{K_{\beta, i}}=\frac{1}{\sqrt[4]{\langle\alpha, \beta\rangle}}\sqrt{K_{\alpha+\beta, i}}$$
and
$$\sqrt{K_{\alpha}}\diamond\sqrt{K_{\beta }}=\frac{1}{\sqrt{\langle\alpha, \beta\rangle}}\sqrt{K_{\alpha+\beta}}.$$
Also the identities in Lemma \ref{lemma hall multiplicatoin of acyclic complexes} are valid for all $\sqrt{K_{\alpha, i}}$ and $\sqrt{K_{\alpha}}$. Obviously one can similarly get a basis of the extended quantum torus $\mathbb{T}_{\Z/t,\frac{1}{2}}^{ac}(\ca)$ as follows.

\begin{proposition}\label{proposition basis of extended torus}
The extended  quantum torus $\mathbb{T}_{\Z/t,\frac{1}{2}}^{ac}(\ca)$ of acyclic $\mathbb{Z}/t$-graded complexes has a basis consisting of the following elements
$$\sqrt{K_{\alpha^0, 0}}\diamond\sqrt{K_{\alpha^1, 1}}\diamond\cdots\diamond\sqrt{K_{\alpha^{t-1}, t-1}},$$
where $\alpha^i\in K_0(\ca)$ for any $i\in\mathbb{Z}/t$.
\end{proposition}

\begin{definition}\label{def extended twisted semi-derived Ringel-Hall algebras}
The extended ($\mathbb{Z}/t$-graded) semi-derived Ringel-Hall algebra $\cs\cd\ch_{\Z/t,\frac{1}{2}}(\ca)$ of $\ca$ is
defined to be the tensor algebra $$\mathbb{T}_{\Z/t,\frac{1}{2}}^{ac}(\ca)\bigotimes\limits_{\mathbb{T}_{\Z/t}^{ac}(\ca)}\cs\cd\ch_{\Z/t}(\ca),$$
with multiplication given by the following rule
$$(\sqrt{[K_1^\bullet]}\diamond\sqrt{[K_1^{'\bullet}]}^{-1}\otimes[M_1^\bullet])\diamond(\sqrt{[K_2^\bullet]}\diamond\sqrt{[K_2^{'\bullet}]}^{-1}\otimes[M_2^\bullet])$$$$
:=\frac{\sqrt{\langle [K_2^{\bullet}], [M_1^\bullet]\rangle}\sqrt{\langle [M_1^\bullet], [K_2^{'\bullet}]\rangle}}{\sqrt{\langle [M_1^\bullet], [K_2^{\bullet}]\rangle}\sqrt{\langle [K_2^{'\bullet}],[M_1^\bullet]\rangle}}(\sqrt{[K_1^\bullet]}\diamond\sqrt{[K_1^{'\bullet}]}^{-1}\diamond\sqrt{[K_2^\bullet]}\diamond\sqrt{[K_2{^{'\bullet}}]}^{-1})\otimes([M_1^\bullet]\diamond[M_2^\bullet]).$$
\end{definition}

So it is easy to see that $\cs\cd\ch_{\Z/t,\frac{1}{2}}(\ca)$ is isomorphic to $\mathbb{T}_{\Z/t,\frac{1}{2}}^{ac}(\ca)\bigotimes\limits_{\mathbb{R}}V$ as linear spaces because of Corollary \ref{cor free module over torus}, where $V$ is the $\mathbb{R}$-linear space generated by the set $\{[\bigoplus\limits_{i=0}^{t-1}U_{A^i,i}]|A^i\in\Iso(\ca),i\in\mathbb{Z}/t\}$.
Similar to the basis of the semi-derived Ringel-Hall algebra $\cs\cd\ch_{\Z/t}(\ca)$ described in Theorem \ref{theorem another basis in SDH(A)}, we have the following proposition.
\begin{proposition}\label{prop basis extend twisted semi-derived}
$\cs\cd\ch_{\Z/t,\frac{1}{2}}(\ca)$ has a basis consisting of the following elements
$$\sqrt{K_{\alpha^0, 0}}\diamond\cdots\diamond\sqrt{K_{\alpha^{t-1}, t-1}}\otimes[U_{A^{t-1},t-1}]\diamond\cdots\diamond[U_{A^{0},0}],$$
where $\alpha^i\in K_0(\ca), [A^i]\in\Iso(\ca)$ for all $i\in\mathbb{Z}/t$.
\end{proposition}

\subsection{Embedding theorem}To prove the main result, we need to show the following isomorphism.
\begin{proposition}\label{proposition Hall number of MHA and DHA}
Let $A^\bullet$ and $B^\bullet$ be $\mathbb{Z}/t$-graded complexes with zero differentials. Then we have the following isomorphism
$$\Ext_{\cc_{\Z/t}(\ca)}^1(A^\bullet, B^\bullet)\cong\Hom_{\cd_{\Z/t}(\ca)}(A^\bullet, B^\bullet[1]).$$
\end{proposition}
\begin{proof}
It is not difficult to see that for each short exact sequence $$\xi^\bullet: 0\rightarrow B^\bullet\xrightarrow{b^\bullet} L^\bullet\xrightarrow{a^\bullet} A^\bullet\rightarrow 0,$$
there exists a homomorphism $$\delta_{\xi^\bullet}=w^\bullet (s^{^\bullet})^{-1}\in\Hom_{\cd_{\Z/t}(\ca)}(A^\bullet, B^\bullet[1])$$
naturally induced by the quasi-isomorphism $s^\bullet=(a^\bullet,0):\Cone(b^\bullet)\rightarrow A^\bullet$ and $w^\bullet=\left(0,1\right): \Cone(b^\bullet)\rightarrow B^\bullet[1]$.
And this induces a homomorphism of abelian groups
$$\Phi: \Ext_{\cc_{\Z/t}(\ca)}^1(A^\bullet, B^\bullet)\rightarrow\Hom_{\cd_{\Z/t}(\ca)}(A^\bullet, B^\bullet[1]),$$
defined by
$$[\xi^\bullet]\mapsto \delta_{\xi^\bullet}$$
for any $[\xi^\bullet]\in\Ext_{\cc_{\Z/t}(\ca)}^1(A^\bullet, B^\bullet)$.

Because the differentials of $A^\bullet$ and $B^\bullet$ are both zero, the differentials of $L^\bullet$ can be described by a homomorphism $c^\bullet\in\Hom_{\cc_{\Z/t}(\ca)}(A^\bullet, B^\bullet[1])$, i.e.,
$$d_{L^\bullet}^i=b^{i+1}c^{i}a^i, \forall i\in\mathbb{Z}/t.$$  Moreover, one should note that all the connecting homomorphisms $\partial^i: A^i\rightarrow B^{i+1}$ of the short exact sequence $\xi^\bullet$ coincide with $c^i$.

For any $[\xi^\bullet]\in\Ext_{\cc_{\Z/t}(\ca)}^1(A^\bullet, B^\bullet)$, if $\delta_{\xi^\bullet}=0$, then there are two homomorphisms $$u^\bullet\in\Hom_{\cd_{\Z/t}(\ca)}(A^\bullet,L^\bullet),~ v^\bullet\in\Hom_{\cd_{\Z/t}(\ca)}(L^\bullet,B^\bullet)$$ such that $a^\bullet u^\bullet=1_{A^\bullet}$ and $v^\bullet b^\bullet=1_{B^\bullet}$ in $\cd_{\Z/t}(\ca)$.  Hence we get the induced isomorphisms
$$H^i(a^\bullet)H^i(u^\bullet)\cong1_{A^i}~\mbox{and}~H^i(v^\bullet)H^i(b^\bullet)\cong1_{B^i}$$
in $\ca$ for each $i\in\mathbb{Z}/t$. On the other hand, it follows from the short exact sequence $\xi^\bullet$ that there is an acyclic $\mathbb{Z}/{3t}$-graded complex of homologies
$$B^0\xrightarrow{H^0(b^\bullet)}H^0(L^\bullet)\xrightarrow{H^0(a^\bullet)}A^0\xrightarrow{c^0}B^1\rightarrow\cdots\rightarrow B^{t-1}\xrightarrow{H^{t-1}(b^\bullet)}H^{t-1}(L^\bullet)\xrightarrow{H^{t-1}(a^\bullet)}A^{t-1}\xrightarrow{c^{t-1}}B^0,$$
so the connecting homomorphism $c^i=0$ for all $i\in\mathbb{Z}/t$, which yields that the differentials $d_{L^\bullet}^i$ of $L^\bullet$ are all zero. In addition, one can get the isomorphism $L^\bullet\cong A^\bullet\oplus B^\bullet$ which means that $\Phi$ is injective.

It remains to show that $\Phi$ is an epimorphism. For any $w^\bullet(s^{^\bullet})^{-1}\in\Hom_{\cd_{\Z/t}(\ca)}(A^\bullet, B^\bullet[1])$, one can assume that it is of the following form
\[\xymatrix{
A^\bullet&X^\bullet\ar@{=>}_{s^\bullet}[l]\ar[r]^{w^\bullet}&B^\bullet[1],
}
\]
where $s^\bullet$ is a quasi-isomorphism. It is easy to see the following short exact sequence
$$0\rightarrow X^\bullet\xrightarrow{\left(\begin{array}{c}s^\bullet\\1\\0\end{array}\right)}A^\bullet\oplus\Cone(1_{X^\bullet})
\xrightarrow{\left(\begin{array}{ccc}1&-s^\bullet&0\\0&0&1\end{array}\right)} \Cone(-s^\bullet)\rightarrow0,$$ so the monomorphism $X^\bullet\xrightarrow{\left(\begin{array}{c}s^\bullet\\1\\0\end{array}\right)}A^\bullet\oplus\Cone(1_{X^\bullet})$  is also a quasi-isomorphism. And then by Proposition \ref{proposition homo-fin acyclic t-periodic complexes} one can obtain the following commuataive diagramm
\[\xymatrix{
0\ar[r]&B^\bullet\ar[r]^{\left(\begin{array}{c}1\\0\end{array}\right)\ \ \ \ \ }\ar@{=}[d]&\Cone(w^\bullet)[-1]\ar@{-->}[d]\ar[r]^{(0~1)}&X^\bullet\ar[r]\ar@{>->}[d]^{\small\left(\begin{array}{c}s^\bullet\\1\\0\end{array}\right)}&0\\
0\ar[r]&B^\bullet\ar@{-->}[r]^{b'^\bullet}&Y^\bullet\ar@{-->}[r]^{a'^\bullet}&A^\bullet\oplus\Cone(1_{X^\bullet})\ar[r]&0
}\]
in $\cc_{\Z/t}(\ca)$, with exact rows. As a result, we have a commutative diagramm
\[\xymatrix{
0\ar[r]&B^\bullet\ar@{-->}[r]^{b^\bullet}\ar@{=}[d]&L^\bullet\ar@{-->}[r]^{a^\bullet}\ar@{-->}[d]&A^\bullet\ar[r]\ar[d]^{\small{\left(\begin{array}{c}1\\0\\0\end{array}\right)}}&0\\
0\ar[r]&B^\bullet\ar[r]^{b'^\bullet}&Y^\bullet\ar[r]^{a'^\bullet\ \ \ \ \ \ }&A^\bullet\oplus\Cone(1_{X^\bullet})\ar[r]&0
}\]
in $\cc_{\Z/t}(\ca)$, induced by the pull-back of $a'^\bullet$ and $\left(\begin{array}{c}1\\0\\0\end{array}\right):A^\bullet\rightarrow A^\bullet\oplus\Cone(1_{X^\bullet})$.

So we get a short exact sequence $\xi^\bullet: 0\rightarrow B^\bullet\xrightarrow{b^\bullet} L^\bullet\xrightarrow{a^\bullet} A^\bullet\rightarrow 0$, and we claim that $\delta_{\xi^\bullet}=w^\bullet(s^{^\bullet})^{-1}$.
In fact, we have the following commutative diagram
\[\xymatrix{
B^\bullet\ar[r]^{\left(\begin{array}{c}1\\0\end{array}\right)\ \ \ \ \ }\ar@{=}[d]&\Cone(w^\bullet)[-1]\ar[d]\ar[r]^{(0~1)}&X^\bullet\ar[r]^{w^\bullet}\ar@{>->}[d]^{\small{\left(\begin{array}{c}s^\bullet\\1\\0\end{array}\right)}}&B^\bullet[1]\ar@{=}[d]\\
B^\bullet\ar[r]^{b'^\bullet}&Y^\bullet\ar[r]^{a'^\bullet}&A^\bullet\oplus\Cone(1_{X^\bullet})\ar[r]&B^\bullet[1]\ar@{=}[d]\\
B^\bullet\ar[r]^{b^\bullet}\ar@{=}[u]&L^\bullet\ar[r]^{a^\bullet}\ar[u]&A^\bullet\ar[r]^{\delta_\xi}\ar@{>->}[u]^{\small{\left(\begin{array}{c}1\\0\\0\end{array}\right)}}&B^\bullet[1],
}\]
in $\cd_{\Z/t}(\ca)$, where the rows are all triangles in $\cd_{\Z/t}(\ca)$. This finishes the proof.
\end{proof}
\begin{theorem}\label{theorem embedding theorem}
There is an embedding $\iota_t:\cd\ch_{\Z/t}(\ca)\rightarrow\cs\cd\ch_{\Z/t,\frac{1}{2}}(\ca)$ for any odd integer $t>0$, defined by
$$Z_A\mapsto\frac{1}{\sqrt{\langle \widehat{A}, \widehat{A}\rangle}}\sqrt{K_{-\widehat{A}}}\otimes U_{A},$$ if $t=1$ and
$$Z_{A}^{[i]}\mapsto\frac{1}{\sqrt[4]{\langle \widehat{A}, \widehat{A}\rangle}}\sqrt{K_{-\widehat{A},i+1}}\diamond\sqrt{K_{\widehat{A},i+2}}\diamond\cdots\diamond\sqrt{K_{\widehat{A},i+t-1}}\diamond\sqrt{K_{-\widehat{A},i+t}}\otimes U_{A,i},$$
if $t>1$, for any $[A]\in\Iso(\ca)$ and $i\in\mathbb{Z}/t$.
\end{theorem}
\begin{proof}
It follows from Proposition \ref{prop basis extend twisted semi-derived} that $\iota_t$ is injective. So it remains to show that $\iota_t$ is a homomorphism of algebras, and we just need to verify it satisfies the relations in Proposition \ref{proposition generators and relations 1-derived hall algebra}, Proposition \ref{proposition generators and relations 3-derived hall algebra} and Proposition \ref{proposition generators and relations t-derived hall algebra} respectively. And we will separate the proof into the following two cases.

Case I. $t=1$. For any $[A], [B]\in\Iso(\ca)$, we have
\begin{eqnarray*}
&&\iota_1(Z_A)\diamond\iota_1(Z_B)\\
&=&\frac{1}{\sqrt{\langle \widehat{A},\widehat{A}\rangle\langle \widehat{B},\widehat{B}\rangle}}(\sqrt{K_{-\widehat{A}}}\otimes U_{A})\diamond(\sqrt{K_{-\widehat{B}}}\otimes U_{B})\\
&=&\frac{1}{\sqrt{\langle \widehat{A},\widehat{A}\rangle\langle \widehat{B},\widehat{B}\rangle}}\frac{\sqrt{\langle \widehat{A},\widehat{B}\rangle}}{\sqrt{\langle \widehat{B},\widehat{A}\rangle}}\sqrt{K_{-\widehat{A}}}\diamond\sqrt{K_{-\widehat{B}}}\otimes U_A\diamond U_B\\
&=&\frac{1}{\sqrt{\langle \widehat{A}+\widehat{B},\widehat{A}+\widehat{B}\rangle}}\sqrt{K_{-(\widehat{A}+\widehat{B})}}\otimes \sqrt{\langle \widehat{A},\widehat{B}\rangle}U_A\diamond U_B\\
&=&\frac{1}{\sqrt{\langle \widehat{A}+\widehat{B},\widehat{A}+\widehat{B}\rangle}}\sqrt{K_{-(\widehat{A}+\widehat{B})}}\otimes \sqrt{\langle \widehat{A},\widehat{B}\rangle}\\
&&\sum\limits_{[C]\in\Iso(\ca)}\left(\sum\limits_{\tiny \begin{array}{c}[L^\bullet]\in\Iso(\cc_{\Z/1}(\ca))\\H(L^\bullet)\cong C\end{array}}\frac{|\Ext^1_{\cc_{\Z/1}(\ca)}(U_A, U_B)_{L^\bullet}|}{|\Hom_{\ca}(A,B)|}\right)\sqrt{\langle \widehat{A}+\widehat{B}-\widehat{C},\widehat{C}\rangle}K_{\frac{1}{2}\left(\widehat{A}+\widehat{B}-\widehat{C}\right)}\diamond U_C
\end{eqnarray*}
By definition one can easily see that
\begin{eqnarray*}
K_{\frac{1}{2}\left(\widehat{A}+\widehat{B}-\widehat{C}\right)}
&=&\sqrt{K_{\frac{1}{2}\left(\widehat{A}+\widehat{B}-\widehat{C}\right)}\oplus K_{\frac{1}{2}\left(\widehat{A}+\widehat{B}-\widehat{C}\right)}}\\
&=&\sqrt{K_{\widehat{A}+\widehat{B}-\widehat{C}}}.
\end{eqnarray*}
So we have that
\begin{eqnarray*}
\iota_1(Z_A)\diamond\iota_1(Z_B)&=&\sum\limits_{[C]\in\Iso(\ca)}\left(\sum\limits_{\tiny \begin{array}{c}[L^\bullet]\in\Iso(\cc_{\Z/1}(\ca))\\H(L^\bullet)\cong C\end{array}}\frac{|\Ext^1_{\cc_{\Z/1}(\ca)}(U_A, U_B)_{L^\bullet}|}{\sqrt{|\Hom_{\ca}(A,B)||\Ext^1_\ca(A,B)|}}\right)\\
&&\frac{\sqrt{\langle \widehat{A}+\widehat{B}-\widehat{C},\widehat{C}\rangle}\sqrt{\langle \widehat{A}+\widehat{B},\widehat{A}+\widehat{B}-\widehat{C}\rangle}}{\sqrt{\langle \widehat{A}+\widehat{B},\widehat{A}+\widehat{B}\rangle}}\sqrt{K_{-\widehat{C}}}\otimes U_C\\
&=&\sum\limits_{[C]\in\Iso(\ca)}\frac{|\Ext^1_{\cd_{\Z/1}(\ca)}(Z_A,Z_B)_C|}{\sqrt{|\Hom_{\cd_{\Z/1}(\ca)}(Z_A,Z_B)|}}\iota_1(Z_C).
\end{eqnarray*}

Case II. $t\geq3$. For any $[A], [B]\in\Iso(\ca)$ and $i\in \mathbb{Z}/t$, we get that
\begin{eqnarray*}
&&\iota_t(Z_{A}^{{[i]}})\diamond\iota_t(Z_{B}^{{[i]}})\\
&=&\frac{1}{\sqrt[4]{\langle \widehat{A},\widehat{A}\rangle\langle \widehat{B},\widehat{B}\rangle}}\left(\sqrt{K_{-\widehat{A},i+1}}\diamond\sqrt{K_{\widehat{A},i+2}}\diamond\cdots\diamond\sqrt{K_{\widehat{A},i-1}}\diamond
\sqrt{K_{-\widehat{A},i}}\otimes U_{A,i}\right)\diamond\\
&&\left(\sqrt{K_{-\widehat{B},i+1}}\diamond\sqrt{K_{\widehat{B},i+2}}\diamond\cdots\diamond\sqrt{K_{\widehat{B},i-1}}\diamond
\sqrt{K_{-\widehat{B},i}}\otimes U_{B,i}\right)\\
&=&\frac{1}{\sqrt[4]{\langle \widehat{A},\widehat{A}\rangle\langle \widehat{B},\widehat{B}\rangle}}\frac{\sqrt{\langle\widehat{A},\widehat{B}}\rangle}{\sqrt{\langle \widehat{B},\widehat{A}\rangle}}\sqrt{K_{-\widehat{A},i+1}}\diamond\sqrt{K_{\widehat{A},i+2}}\diamond\cdots\diamond\sqrt{K_{\widehat{A},i-1}}\diamond
\sqrt{K_{-\widehat{A},i}}\diamond\\
&&\sqrt{K_{-\widehat{B},i+1}}\diamond\sqrt{K_{\widehat{B},i+2}}\diamond\cdots\diamond\sqrt{K_{\widehat{B},i-1}}\diamond
\sqrt{K_{-\widehat{B},i}}\otimes U_{A,i}\diamond U_{B,i}
\end{eqnarray*}
\begin{eqnarray*}&=&\frac{\sqrt{\langle\widehat{A},\widehat{B}}\rangle}{\sqrt[4]{\langle \widehat{A},\widehat{A}\rangle\langle \widehat{B},\widehat{B}\rangle\langle \widehat{B},\widehat{A}\rangle}}\sqrt{K_{-(\widehat{A}+\widehat{B}),i+1}}\diamond\sqrt{K_{\widehat{A},i+2}}\diamond\cdots\diamond\sqrt{K_{\widehat{A},i-1}}\diamond
\sqrt{K_{-\widehat{A},i}}\diamond\\
&&\sqrt{K_{\widehat{B},i+2}}\diamond\cdots\diamond\sqrt{K_{\widehat{B},i-1}}\diamond
\sqrt{K_{-\widehat{B},i}}\otimes U_{A,i}\diamond U_{B,i}\\
&=&\frac{\sqrt{\langle\widehat{A},\widehat{B}}\rangle}{\sqrt[4]{\langle \widehat{A},\widehat{A}\rangle\langle \widehat{B},\widehat{B}\rangle\langle \widehat{B},\widehat{A}\rangle}}\sqrt{K_{-(\widehat{A}+\widehat{B}),i+1}}\diamond\sqrt{K_{\widehat{A}+\widehat{B},i+2}}\diamond\cdots\diamond
\sqrt{K_{-(\widehat{A}+\widehat{B}),i-2}}\diamond\\
&&\sqrt{K_{\widehat{A},i-1}}\diamond
\sqrt{K_{-\widehat{A},i}}\diamond\sqrt{K_{\widehat{B},i-1}}\diamond
\sqrt{K_{-\widehat{B},i}}\otimes U_{A,i}\diamond U_{B,i}\\
&=&\frac{\sqrt{\langle\widehat{A},\widehat{B}}\rangle}{\sqrt[4]{\langle \widehat{A}+\widehat{B},\widehat{A}+\widehat{B}\rangle}}\sqrt{K_{-(\widehat{A}+\widehat{B}),i+1}}\diamond\cdots\diamond\sqrt{K_{-(\widehat{A}+\widehat{B})_,i}}\otimes U_{A,i}\diamond U_{B,i}\\
&=&\sum\limits_{[C]\in\Iso(\ca)}\frac{|\Ext_{\ca}^{1}(A,B)_{C}|}{\sqrt{|\Hom_{\ca}(A,B)||\Ext^1_\ca(A,B)|}}\iota_t(Z_{C}^{{[i]}}),
\end{eqnarray*}
and
\begin{eqnarray*}
&&\iota_t(Z_{B}^{{[i]}})\diamond\iota_t(Z_{A}^{{[i+1]}})\\
&=&\frac{1}{\sqrt[4]{\langle \widehat{B},\widehat{B}\rangle\langle \widehat{A},\widehat{A}\rangle}}\left(\sqrt{K_{-\widehat{B},i+1}}\diamond\sqrt{K_{\widehat{B},i+2}}\diamond\cdots\diamond\sqrt{K_{\widehat{B},i-1}}\diamond
\sqrt{K_{-\widehat{B},i}}\otimes U_{B,i}\right)\diamond\\
&&\left(\sqrt{K_{-\widehat{A},i+2}}\diamond\sqrt{K_{\widehat{A},i+3}}\diamond\cdots\diamond\sqrt{K_{\widehat{A},i}}\diamond
\sqrt{K_{-\widehat{A},i+1}}\otimes U_{A,i+1}\right)\\
&=&\frac{1}{\sqrt[4]{\langle \widehat{B},\widehat{B}\rangle\langle \widehat{A},\widehat{A}\rangle}}\frac{1}{\sqrt{( \widehat{A},\widehat{B})}}\sqrt{K_{-\widehat{B},i+1}}\diamond\sqrt{K_{\widehat{B},i+2}}\diamond\cdots\diamond\sqrt{K_{\widehat{B},i-1}}\diamond
\sqrt{K_{-\widehat{B},i}}\diamond\\
&&\sqrt{K_{-\widehat{A},i+2}}\diamond\sqrt{K_{\widehat{A},i+3}}\diamond\cdots\diamond\sqrt{K_{\widehat{A},i}}\diamond
\sqrt{K_{-\widehat{A},i+1}}\otimes U_{B,i}\diamond U_{A,i+1}\\
&=&\frac{1}{\sqrt[4]{\langle \widehat{B},\widehat{B}\rangle\langle \widehat{A},\widehat{A}\rangle}}\frac{1}{\sqrt{( \widehat{A},\widehat{B})}}\frac{1}{\sqrt{\langle\widehat{B},\widehat{B}\rangle}}\sqrt{K_{\widehat{B},i+2}}\diamond\cdots\diamond\sqrt{K_{\widehat{B},i-1}}\diamond
\sqrt{K_{-\widehat{B},i}}\diamond\sqrt{K_{-\widehat{B},i+1}}\diamond\\
&&\sqrt{K_{-\widehat{A},i+2}}\diamond\sqrt{K_{\widehat{A},i+3}}\diamond\cdots\diamond\sqrt{K_{\widehat{A},i}}\diamond
\sqrt{K_{-\widehat{A},i+1}}\otimes U_{B,i}\diamond U_{A,i+1}\\
&=&\frac{1}{\sqrt[4]{\langle \widehat{B},\widehat{B}\rangle\langle \widehat{A},\widehat{A}\rangle}}\frac{1}{\sqrt{( \widehat{A},\widehat{B})}}\frac{\sqrt[4]{\langle\widehat{A},\widehat{B}}\rangle}{\sqrt{\langle\widehat{B},\widehat{B}\rangle}}
\sqrt{K_{\widehat{B}-\widehat{A},i+2}}\diamond\sqrt{K_{-\widehat{B},i+3}}\diamond\cdots\diamond
\sqrt{K_{-\widehat{B},i}}\diamond\sqrt{K_{-\widehat{B},i+1}}\diamond\\
&&\sqrt{K_{\widehat{A},i+3}}\diamond\cdots\diamond\sqrt{K_{\widehat{A},i}}\diamond
\sqrt{K_{-\widehat{A},i+1}}\otimes U_{B,i}\diamond U_{A,i+1}\\
&=&\frac{1}{\sqrt[4]{\langle \widehat{B},\widehat{B}\rangle\langle \widehat{A},\widehat{A}\rangle}}\frac{1}{\sqrt{( \widehat{A},\widehat{B})}}\frac{\sqrt[4]{\langle\widehat{A},\widehat{B}}\rangle}{\sqrt{\langle\widehat{B},\widehat{B}\rangle}}
\sqrt{K_{\widehat{B}-\widehat{A},i+2}}\diamond\sqrt{K_{-(\widehat{B}-\widehat{A}),i+3}}\diamond\cdots\diamond\sqrt{K_{\widehat{B}-\widehat{A},i-1}}\diamond\\
&&\sqrt{K_{-\widehat{B},i}}\diamond\sqrt{K_{-\widehat{B},i+1}}\diamond\sqrt{K_{\widehat{A},i}}\diamond\sqrt{K_{-\widehat{A},i+1}}\otimes U_{B,i}\diamond U_{A,i+1}\\
&=&\frac{\sqrt[4]{\langle\widehat{A},\widehat{B}}\rangle\sqrt[4]{\langle\widehat{B},\widehat{A}}\rangle}{\sqrt[4]{\langle \widehat{B},\widehat{B}\rangle\langle \widehat{A},\widehat{A}\rangle}\sqrt{(\widehat{A},\widehat{B})}\sqrt{\langle \widehat{B},\widehat{B}\rangle}}
\sqrt{K_{(\widehat{B}-\widehat{A}),i+2}}\diamond\sqrt{K_{(\widehat{A}-\widehat{B}),i+3}}\diamond\cdots\diamond\sqrt{K_{(\widehat{A}-\widehat{B}),i}}\diamond
\end{eqnarray*}
\begin{eqnarray*}&&\sqrt{K_{-(\widehat{A}+\widehat{B}),i+1}}\otimes(U_{B,i}\diamond U_{A,i+1})\\
&=&\sum\limits_{[M],[N]\in\Iso(\ca)}\gamma_{AB}^{MN}\frac{a_Aa_B}{a_Ma_N}\frac{\langle\widehat{B}-\widehat{M},\widehat{M}\rangle}{\sqrt[4]{\langle \widehat{A}+\widehat{B}, \widehat{A}+\widehat{B}\rangle}\sqrt{\langle \widehat{B},\widehat{B}\rangle}}\sqrt{K_{(\widehat{B}-\widehat{A}),i+2}}\diamond\cdots\diamond\sqrt{K_{-(\widehat{A}+\widehat{B}),i+1}}\diamond \\
&&K_{\widehat{B}-\widehat{M},i+1}\otimes U_{N,i+1}\diamond U_{M,i}\\
&=&\sum\limits_{[M],[N]\in\Iso(\ca)}\gamma_{AB}^{MN}\frac{a_Aa_B}{a_Ma_N}\frac{\langle\widehat{B}-\widehat{M},\widehat{M}\rangle\sqrt{\langle \widehat{A}+\widehat{B},\widehat{B}-\widehat{M}\rangle}}{\sqrt[4]{\langle \widehat{A}+\widehat{B}, \widehat{A}+\widehat{B}\rangle}\sqrt{\langle \widehat{B},\widehat{B}\rangle}}\sqrt{K_{\widehat{M}-\widehat{N},i+2}}\diamond\cdots\diamond\sqrt{K_{-(\widehat{M}-\widehat{N}),i}}\diamond
\\
&&\sqrt{K_{-(\widehat{N}+\widehat{M}),i+1}}\otimes U_{N,i+1}\diamond U_{M,i}.
\end{eqnarray*}
In addition, we obtain that
\begin{eqnarray*}
\frac{\langle\widehat{B}-\widehat{M},\widehat{M}\rangle\sqrt{\langle \widehat{A}+\widehat{B},\widehat{B}-\widehat{M}\rangle}}{\sqrt[4]{\langle \widehat{A}+\widehat{B}, \widehat{A}+\widehat{B}\rangle}\sqrt{\langle \widehat{B},\widehat{B}\rangle}}&=&\frac{\langle\widehat{B}-\widehat{M},\widehat{M}\rangle}{\sqrt{\langle \widehat{B},\widehat{B}\rangle}\sqrt[4]{\langle \widehat{A}+\widehat{B}, \widehat{A}+\widehat{B}-2\widehat{B}+2\widehat{M}\rangle}}
\\
&=&\frac{\langle\widehat{B}-\widehat{M},\widehat{M}\rangle}{\sqrt{\langle \widehat{B},\widehat{B}\rangle}\sqrt[4]{\langle \widehat{A}+\widehat{B}, \widehat{M}+\widehat{N}\rangle}}\\
&=&\frac{\sqrt{\langle\widehat{B},\widehat{M}+\widehat{N}-\widehat{A}\rangle}}{\sqrt[4]{\langle \widehat{A}+\widehat{B}, \widehat{M}+\widehat{N}\rangle}\langle\widehat{M},\widehat{M}\rangle}\\
&=&\frac{1}{\sqrt[4]{\langle \widehat{A}-\widehat{B},\widehat{M}+\widehat{N}\rangle}\sqrt{\langle \widehat{B},\widehat{A}\rangle}\langle\widehat{M},\widehat{M}\rangle}\\
&=&\frac{1}{\sqrt[4]{\langle \widehat{N}-\widehat{M},\widehat{M}+\widehat{N}\rangle}\sqrt{\langle\widehat{M},\widehat{M}\rangle\langle\widehat{M},\widehat{N}\rangle}}
\frac{\sqrt{\langle\widehat{M},\widehat{N}\rangle}}{\sqrt{\langle \widehat{B},\widehat{A}\rangle\langle\widehat{M},\widehat{M}\rangle}}\\
&=&\frac{1}{\sqrt[4]{\langle \widehat{N}+\widehat{M},\widehat{M}+\widehat{N}\rangle}}\frac{\sqrt{\langle\widehat{M},\widehat{N}\rangle}}{\sqrt{\langle \widehat{B},\widehat{A}\rangle\langle\widehat{M},\widehat{M}\rangle}}
\end{eqnarray*}
On the other hand, we get that
\begin{eqnarray*}
&&\sum\limits_{[M],[N]\in\Iso(\ca)}\gamma_{AB}^{MN}\frac{a_Aa_B}{a_Ma_N}\frac{1}{\sqrt{\langle \widehat{B}, \widehat{A}\rangle\langle\widehat{N},\widehat{M}\rangle} } \iota_t(Z_{N}^{{[i+1]}})\diamond\iota_t(Z_{M}^{{[i]}})\\
&=&\sum\limits_{[M],[N]\in\Iso(\ca)}\gamma_{AB}^{MN}\frac{a_Aa_B}{a_Ma_N}\frac{1}{\sqrt{\langle \widehat{B}, \widehat{A}\rangle\langle\widehat{N},\widehat{M}\rangle}\sqrt[4]{\langle \widehat{N},\widehat{N}\rangle\langle \widehat{M},\widehat{M}\rangle}}\\
&&(\sqrt{K_{-\widehat{N},i+2}}\diamond\sqrt{K_{\widehat{N},i+3}}\diamond\cdots\diamond\sqrt{K_{\widehat{N},i}}\diamond
\sqrt{K_{-\widehat{N},i+1}}\otimes U_{N,i+1})\diamond\\
&&(\sqrt{K_{-\widehat{M},i+1}}\diamond\sqrt{K_{\widehat{M},i+2}}\diamond\cdots\diamond\sqrt{K_{\widehat{M},i-1}}\diamond
\sqrt{K_{-\widehat{M},i}}\otimes U_{M,i})\\
&=&\sum\limits_{[M],[N]\in\Iso(\ca)}\gamma_{AB}^{MN}\frac{a_Aa_B}{a_Ma_N}\frac{\sqrt{\langle \widehat{N},\widehat{M}\rangle\langle \widehat{M},\widehat{N}\rangle}}{\sqrt{\langle \widehat{B}, \widehat{A}\rangle\langle\widehat{N},\widehat{M}\rangle}\sqrt[4]{\langle \widehat{N},\widehat{N}\rangle\langle \widehat{M},\widehat{M}\rangle}}\\
\end{eqnarray*}
\begin{eqnarray*}&&\sqrt{K_{-\widehat{N},i+2}}\diamond\sqrt{K_{\widehat{N},i+3}}\diamond\cdots\diamond\sqrt{K_{\widehat{N},i}}\diamond
\sqrt{K_{-\widehat{N},i+1}}\diamond\\
&&\sqrt{K_{-\widehat{M},i+1}}\diamond\sqrt{K_{\widehat{M},i+2}}\diamond\cdots\diamond\sqrt{K_{\widehat{M},i-1}}\diamond
\sqrt{K_{-\widehat{M},i}}\otimes  U_{N,i+1}\diamond U_{M,i}\\
&=&\sum\limits_{[M],[N]\in\Iso(\ca)}\gamma_{AB}^{MN}\frac{a_Aa_B}{a_Ma_N}\frac{\sqrt{\langle \widehat{M},\widehat{N}\rangle}}{\sqrt{\langle \widehat{B}, \widehat{A}\rangle}\sqrt[4]{\langle \widehat{N},\widehat{N}\rangle\langle \widehat{M},\widehat{M}\rangle\langle\widehat{N},\widehat{M}\rangle}}\sqrt{K_{-\widehat{N},i+2}}
\diamond\sqrt{K_{\widehat{N},i+3}}\diamond\cdots\diamond\\
&&\sqrt{K_{\widehat{N},i}}\diamond\sqrt{K_{-(\widehat{M}+\widehat{N}),i+1})}\diamond\sqrt{K_{\widehat{M},i+2}}\diamond\cdots\diamond\sqrt{K_{\widehat{M},i-1}}\diamond
\sqrt{K_{-\widehat{M},i}}\otimes U_{N,i+1}\diamond U_{M,i}\\
&=&\sum\limits_{[M],[N]\in\Iso(\ca)}\gamma_{AB}^{MN}\frac{a_Aa_B}{a_Ma_N}\frac{\sqrt{\langle \widehat{M},\widehat{N}\rangle}}{\sqrt{\langle \widehat{B}, \widehat{A}\rangle}\sqrt[4]{\langle \widehat{N},\widehat{N}\rangle\langle \widehat{M},\widehat{M}\rangle\langle\widehat{N},\widehat{M}\rangle\langle \widehat{M},\widehat{N}+\widehat{M}\rangle}}\sqrt{K_{\widehat{M}-\widehat{N},i+2}}\diamond\cdots\diamond\\
&&\sqrt{K_{-\widehat{N},i}}\diamond\sqrt{K_{-(\widehat{M}+\widehat{N}),i+1})}\diamond\sqrt{K_{-\widehat{M},i+3}}\diamond\cdots\diamond\sqrt{K_{-\widehat{M},i}}\otimes U_{N,i+1}\diamond U_{M,i}
\\
&=&\sum\limits_{[M],[N]\in\Iso(\ca)}\gamma_{AB}^{MN}\frac{a_Aa_B}{a_Ma_N}\frac{\sqrt{\langle \widehat{M},\widehat{N}\rangle}}{\sqrt{\langle \widehat{B}, \widehat{A}\rangle}\sqrt[4]{\langle \widehat{N},\widehat{N}\rangle\langle \widehat{M},\widehat{M}\rangle\langle\widehat{N},\widehat{M}\rangle\langle \widehat{M},\widehat{N}+\widehat{M}\rangle}}\\
&&\sqrt{K_{\widehat{M}-\widehat{N},i+2}}\diamond\sqrt{K_{\widehat{N}-\widehat{M},i+3}}\diamond\cdots\diamond
\sqrt{K_{\widehat{M}-\widehat{N},i-1}}\diamond\sqrt{K_{\widehat{N},i}}\diamond\sqrt{K_{-(\widehat{M}+\widehat{N}),i+1}}\diamond\sqrt{K_{-\widehat{M},i}}\otimes\\
&&U_{N,i+1}\diamond U_{M,i}\\
&=&\sum\limits_{[M],[N]\in\Iso(\ca)}\gamma_{AB}^{MN}\frac{a_Aa_B}{a_Ma_N}\frac{\sqrt{\langle \widehat{M},\widehat{N}\rangle}}{\sqrt{\langle \widehat{B}, \widehat{A}\rangle}\sqrt[4]{\langle \widehat{N},\widehat{N}\rangle\langle \widehat{M},\widehat{M}\rangle\langle\widehat{N},\widehat{M}\rangle\langle \widehat{M},\widehat{N}+\widehat{M}\rangle}}\frac{\sqrt[4]{\langle \widehat{N},\widehat{M}\rangle}}{\sqrt[4]{\langle\widehat{N}+\widehat{M},\widehat{M}}\rangle}\\
&&\sqrt{K_{\widehat{M}-\widehat{N},i+2}}\diamond\cdots\diamond\sqrt{K_{-(\widehat{M}-\widehat{N}),i}}\diamond
\sqrt{K_{-(\widehat{N}+\widehat{M}),i+1}}\otimes U_{N,i+1}\diamond U_{M,i}\\
&=&\sum\limits_{[M],[N]\in\Iso(\ca)}\gamma_{AB}^{MN}\frac{a_Aa_B}{a_Ma_N}\frac{1}{\sqrt[4]{\langle \widehat{N}+\widehat{M},\widehat{M}+\widehat{N}\rangle}}\frac{\sqrt{\langle\widehat{M},\widehat{N}\rangle}}{\sqrt{\langle \widehat{B},\widehat{A}\rangle\langle\widehat{M},\widehat{M}\rangle}}\\
&&\sqrt{K_{\widehat{M}-\widehat{N},i+2}}\diamond\cdots\diamond\sqrt{K_{-(\widehat{M}-\widehat{N}),i}}\diamond
\sqrt{K_{-(\widehat{N}+\widehat{M}),i+1}}\otimes U_{N,i+1}\diamond U_{M,i}.
\end{eqnarray*}
So $\iota_t$ satisfies the relations (\ref{prop relation 1 in DH3}) and (\ref{prop relation 2 in DH3}).

For the relation (\ref{relation 3 in DHt}), we have
\begin{eqnarray*}
&&\iota_t(Z_A^{[i]})\diamond\iota_t(Z_B^{[j]})\\
&=&\frac{1}{\sqrt[4]{\langle \widehat{A},\widehat{A}\rangle\langle \widehat{B},\widehat{B}\rangle}}\left(\sqrt{K_{-\widehat{A},i+1}}\diamond\sqrt{K_{\widehat{A},i+2}}\diamond\cdots\diamond\sqrt{K_{\widehat{A},i-1}}\diamond
\sqrt{K_{-\widehat{A},i}}\otimes U_{A,i}\right)\diamond\\
&&\left(\sqrt{K_{-\widehat{B},j+1}}\diamond\sqrt{K_{\widehat{B},j+2}}\diamond\cdots\diamond\sqrt{K_{\widehat{B},j-1}}\diamond
\sqrt{K_{-\widehat{B},j}}\otimes U_{B,j}\right)\\
&=&\frac{\sqrt{\langle \widehat{B},\widehat{A}\rangle^{(-1)^{j-i}}\langle \widehat{A},\widehat{B}\rangle^{(-1)^{j-i}}}}{\sqrt[4]{\langle \widehat{A},\widehat{A}\rangle\langle \widehat{B},\widehat{B}\rangle}}
\left(\sqrt{K_{-\widehat{B},j+1}}\diamond\sqrt{K_{\widehat{B},j+2}}\diamond\cdots\diamond\sqrt{K_{\widehat{B},j-1}}\diamond
\sqrt{K_{-\widehat{B},j}}\otimes U_{B,j}\right)\diamond\\
&&\left(\sqrt{K_{-\widehat{A},i+1}}\diamond\sqrt{K_{\widehat{A},i+2}}\diamond\cdots\diamond\sqrt{K_{\widehat{A},i-1}}\diamond
\sqrt{K_{-\widehat{A},i}}\otimes U_{A,i}\right)\\
&=&\sqrt{(\widehat{B},\widehat{A})^{(-1)^{j-i}}}\iota_t(Z_B^{[j]})\diamond\iota_t(Z_A^{[i]}).
\end{eqnarray*}
This completes the proof.
\end{proof}

\end{document}